
\magnification=1200
\documentstyle{amsppt}
\pageheight{9truein} \pagewidth{6.5truein} \pageno=1
        
\def\yawn #1{\nobreak\par\nobreak
\vskip #1 plus .1 in \penalty 200}

\redefine\le{\leqslant} \redefine\ge{\geqslant}
\def\op{\operatorname}
\def\Z{\Bbb Z} \def\Q{{\Bbb Q}} \def\R{\Bbb R}
\def\C{\Bbb C} \redefine\P{\Bbb P}

\def\cal{\Cal}
\redefine\phi{\varphi} 
\def\bs{\backslash}

\def\SL{\operatorname{SL}}
\def\GL{\operatorname{GL}}
\def\geod{\medspace \text{\bf o} \medspace}

\def\qT{{{}_\Q}{\bold T}}
\def\bbs #1{\bold b^{(#1)}}
\def\rs #1{r^{(#1)}}
\def\rts #1{r_t^{(#1)}}
\def\Ms #1{M^{(#1)}}
\def\F{{\Cal F}}

\hyphenation{con-figu-ration}
\hyphenation{con-figu-rations}

\NoBlackBoxes


\topmatter
\title Cohomology at Infinity \\ And the Well-Rounded Retract \\ For
General Linear Groups
\endtitle
\rightheadtext{Cohomology at Infinity and the Well-Rounded Retract}
\author Avner Ash and Mark McConnell \endauthor
\address Dept.~of Mathematics, Ohio State University, Columbus, OH
43210-1174 \endaddress
\email ash\@math.ohio-state.edu \endemail
\address Oklahoma State University, Math, 401 Math Sciences, Stillwater,
OK 74078-1058 \endaddress
\email mmcconn\@math.okstate.edu \endemail
\thanks
Research partially supported by NSA grant MDA-904-94-H-2030 (first author)
and NSF grants DMS-9206659 and 9401460 (second author).  This
manuscript is submitted for publication with the understanding that
the United States government is authorized to reproduce and distribute
reprints.\endthanks
\endtopmatter

\document


\heading Introduction \endheading

\subheading{(0.1)} Let $\bold G$ be a reductive algebraic group
defined over $\Q$, and let $\Gamma$ be an arithmetic subgroup of
$\bold G(\Q)$.  Let $X$ be the symmetric space for $\bold G(\R)$, and
assume $X$ is contractible.  Then the cohomology (mod torsion) of the
space $X/\Gamma$ is the same as the cohomology of $\Gamma$.  In turn,
$X/\Gamma$ will have the same cohomology as $W/\Gamma$, if $W$ is a
``spine'' in $X$.  This means that $W$ (if it exists) is a deformation
retract of $X$ by a $\Gamma$-equivariant deformation retraction, that
$W/\Gamma$ is compact, and that $\dim W$ equals the virtual
cohomological dimension (vcd) of $\Gamma$.  Then $W$ can be given the
structure of a cell complex on which $\Gamma$ acts cellularly, and the
cohomology of $W/\Gamma$ can be found combinatorially.

Spines have been found for many groups $\bold G$ (see~(2.6) below, and
\cite{A1} \cite{M-M0} \cite{M-M1}).  This paper concerns the case
where $\bold G$ is the restriction of scalars of a general linear
group over a number field $k$ with ring of integers $\cal O$.  This
means $\bold G(\Q) = \GL_n(k)$.  We shall take $\Gamma$ to be a
subgroup of finite index in $\GL_n(\cal O)$, or more generally in
$\GL(\cal P)$, where $\cal P$ is a projective $\cal O$-module of
rank~$n$.  In~\cite{A3}, Ash found a spine $W$ for these $\bold G$,
calling $W$ the {\sl well-rounded retract}.  (This generalized
\cite{Sou2}.) The retract has been used in computations:
see~\cite{Sou1} \cite{A-G-G} \cite{A-M1} \cite{A-M2} \cite{vG-T}.  We
remark that $k=\Q$, $\Gamma \subseteq \GL_n(\Z)$ for $n = 2,3,4$ still
provide our main cases of computational interest.  For imaginary
quadratic fields~$k$ and $n=2$, see~\cite{Men} \cite{S-V} \cite{V}.
For real quadratic~$k$ and $n=2$, see~\cite{B}.

\subheading{(0.2)} The present paper extends \cite{A3} to deal with
the ``cusps'' of $X/\Gamma$, in a way which we now describe.

Borel and Serre \cite{B-S} introduced a bordification $\bar X$ of $X$
such that $\bar X/\Gamma$ is a compactification of $X/\Gamma$.  The
space $\bar X/\Gamma$ has the same homotopy type as $X/\Gamma$, and
therefore has the same homotopy type as $W/\Gamma$.  The boundary of
$\bar X /\Gamma$ is a union of finitely many ``faces'', one for each
equivalence class mod~$\Gamma$ of parabolic $\Q$-subgroups of~$\bold G$.

The Borel-Serre boundary gives a geometric approach for distinguishing
the cuspidal cohomology within the whole cohomology of $\Gamma$.  A
cohomology class on $\bar X/\Gamma$ with coefficients in~$\C$ can be
viewed, via the de Rham theorem, as coming from a $\Gamma$-invariant
automorphic form on $\bold G(\R)$.  If the class has non-zero
restriction to the boundary, it is certainly not a cusp form.  The
converse is not completely true, but the difference between the
``interior cohomology'' (the span of classes that restrict to zero on
the boundary) and the ``cuspidal cohomology'' (the span of classes
whose associated automorphic forms are cuspidal) is fairly well
understood, thanks to the theory of Eisenstein series.  In many
practical examples, like the $\GL_4(\Z)$ example mentioned below,
there is no difference; to compute the cuspidal cohomology in certain
degrees, we simply compute the restriction of the cohomology to the
boundary.

A related question is the existence of ``ghost classes''---those
classes which restrict to zero on each face of the boundary, but do not
restrict to zero on the whole boundary.  Harder has given examples of
ghost classes \cite{H}; it is of interest to know how widespread the
phenomenon is.

In specific cases, then, one may need to compute the cohomology of the
boundary of $\bar X/\Gamma$, the cohomology of individual
boundary faces, and the restriction map from $H^*(\bar
X/\Gamma)$ to these groups.

This paper builds a technology for these computations, when $\bold G$
is the general linear group in~(0.1) and $W$ is the well-rounded
retract of \cite{A3}.  To each Borel-Serre boundary face in $\bar X$,
we associate a closed subcomplex $W'$ of $W$.  Let $\Gamma'$ be
the stabilizer of $W'$ in $\Gamma$.  We prove that the inclusion
$W'/\Gamma' \hookrightarrow W/\Gamma$ induces the same map on
cohomology as the inclusion of the corresponding boundary face into
$\bar X/\Gamma$.  Our main result, Theorem~9.5, is that there is a
spectral sequence converging to $H^*(\partial \bar X/\Gamma)$ whose
$E_1$ term is the direct sum of the cohomology of the various
$W'/\Gamma'$.  There are only finitely many terms in this direct sum.
There is a map of spectral sequences, starting with a natural map from
$H^*(W/\Gamma)$ to the $E_1$ term, that induces the canonical
restriction map $H^*(\bar X/\Gamma) \to H^*(
\partial \bar
X/\Gamma)$ up to canonical isomorphism.  There is a dual statement in
homology~(10.2).  We do not assume that $\Gamma$ is torsion-free.  We
use (co)homology with coefficients in any fixed abelian group.

We can also find the maps on (co)homology induced by the inclusion of
a single boundary face into $\bar X/\Gamma$, without any spectral
sequences~(10.3).

Each $W'/\Gamma'$ is a finite cell complex, so its cohomology can be
worked out combinatorially.  See~(6.3) for an example.  We emphasize
that each $W'$ is a subcomplex of $W$, so that we can work with it in
terms already available once we have computed $W$ in any given case.

In a future paper, we will use this method to study the congruence
subgroups $\Gamma_0(p)$ of prime level~$p$ in $\GL_4(\Z)$.  For a
range of~$p$, we will compute $H^5(\Gamma_0(p), \C)$, and will isolate
the subspace of classes that restrict to zero on the boundary.  We
will compute the Hecke eigenvalues of these classes, using either the
new techniques of Paul Gunnells, or the new techniques of
R.~MacPherson and the second author.  We will try to determine whether
any of the cusp forms we find are functorial lifts from smaller groups
such as $\op{Sp}_4$ or $\op{O}_4$.

\subheading{(0.3)} We have long wished to see spelled out in detail
the relationship between the well-rounded retract and the Borel-Serre
compactification.  The retraction $r:X\to W$ carries any $x\in X$ into
$W$ by a combination of geodesic actions depending on~$x$---see~(4.2).
Each Borel-Serre boundary component $e$ in $\bar X$ has a tubular
neighborhood $\bar N$, homeomorphic to $e\times [0,1)^l$ for
appropriate~$l$, where the fibers $[0,1)^l$ are given by the geodesic
action.  The retraction~$r$ is constant on the interior of each fiber,
and carries the interior of $\bar N$ onto the corresponding cell
complex $W'$.

The main idea of the paper is to extend the retraction $r: X\to W$ to
a continuous $\Gamma$-equivariant map $\bar r:\bar X\to W$.  It is
relatively easy to see that the continuous extension exists and is
$\Gamma$-equivariant.
However, $r = r_1$ really comes from a deformation retraction $r_t :
X\times [0,1] \to X$, and it is harder to prove that $\bar r$ comes from a
deformation retraction of $\bar X$ onto $W$.  A simple example will
illustrate the problem.  The map $f_t : \R\times [0,1]\to \R$ given by
$(x,t)\mapsto (1-t)\cdot x$ is a deformation retraction of $\R$ onto
$\{0\}$.  The compactification $[-\infty,\infty]$ of $\R$ certainly
has a deformation retraction onto $\{0\}$.  But $f_t$ does not extend
to the deformation retraction for $[-\infty,\infty]$.

To fix the problem with $[-\infty,\infty]$, we can first use any
deformation retraction onto $[-x_0, x_0]$ for some large $x_0$, and
then use the old $f_t$ to retract $[-x_0, x_0]$ to $\{0\}$.  We do
something similar to prove $\bar r$ is a deformation retraction.  We
use Saper's recent construction~\cite{S} of a ``central tile'' $X_0
\subset X$ and a deformation retraction $r_{t,S}$ of $\bar X$ onto
$X_0$.  This $X_0$ is a large codimension-zero subset of $X$, which we
may take arbitrarily close to $\partial \bar X$ in all
directions.  We show that composing our $r_t$ with Saper's $r_{t,S}$
gives a deformation retraction $R_t:\bar X\times[0,1] \to \bar X$ with
$\bar r = R_1$.

\comment
\footnote{To the tune of The Hokey Pokey: ``You do a little
Saper/ And you turn your $L$'s around/ That's what it's all about.''}
\endcomment

By what we have said, $\bar r$ gives a homotopy equivalence between $\bar
X/\Gamma$ and $W/\Gamma$.  We show it also gives a homotopy
equivalence between any Borel-Serre boundary component $\bar e/\Gamma'$ and
the corresponding cell complex $W'/\Gamma'$.  It is then routine to
show that $\bar r$ gives an isomorphism between a standard spectral
sequence for $H^*(\partial \bar X/\Gamma)$ and the spectral sequence
described in~(0.2).  This proves our main theorem.

The sections of the paper follow this outline.  In Sections~1--4 we
summarize basic material on algebraic groups, lattices, the
well-rounded retract and retraction, and the Borel-Serre
compactification.  None of this material is new, but we must recall it
to establish our notation, and because we emphasize some less familiar
aspects of it---marked lattices, the flag of successive minima, the
orthogonal scaling group (which is the geodesic action from a
lattice's point of view), and the connections between $r_t$ and the
geodesic action.  Section~5 summarizes our notation.  In Section~6 we
define the sets $W_{\cal F}$ (called $W'$ in this introduction).  In
Section~7 we define the tubular neighborhoods $\bar N_{\cal F}$
(called $\bar N$ here); this is the technical heart of the paper,
where we explore the connections between the well-rounded retraction
and the geodesic action.  Section~8 summarizes the material
from~\cite{S} that we need, and contains our results on $R_t$ and
$\bar r$.  In Section~9, we use $\bar r$ to prove our main theorem on
spectral sequences.

\subheading{(0.4)} The results in this paper could be extended in
several ways.  First, \cite{A3} defines $W$ for $\GL_n$ of any finite
dimensional division algebra over $\Q$, not necessarily commutative.
Our results could be generalized to that case, presumably with only
technical changes.  Second, our arguments work for $\SL$ as well as
$\GL$, with minor changes~(10.5).  See~(10.4) for
comments on $\Gamma$-equivariant cohomology, and~(10.6)--(10.7) for
other connections with reduction theory.

Parts (1.6), (2.6), and (6.3) present the case $k=\Q$ in more
detail.  We hope they help guide the reader through the technicalities
that involve structure over $k$ versus structure over $\Q$.

We are grateful to J. R. Myers and Steve Zucker for helpful
conversations.  We especially thank Chun-Nip Lee and Leslie Saper.  We
are grateful to the referee for several comments.  We thank Bob
MacPherson and Dan Grayson for introducing us to some of these ideas,
and Connie Chen for help with the illustration.  We thank the NSA and
the NSF for support.


\heading Section 1---Background on Algebraic Groups \endheading

\subheading{(1.1)} We fix a number field $k$ of
degree $d$ over $\Q$, with ring of integers~$\cal O$.  Let $v$ run
through the archimedean places of $k$; the completion $k_v$ is either
$\R$ or $\C$.  Let $\R_+$ be the positive real numbers.

\subheading{(1.2)} Fix an integer $n\ge 1$.  The algebraic group
$\GL_n$ is defined over $\Q$, and hence over~$k$.  Let $\bold G$ be
the restriction of scalars of $\GL_n$ from $k$ down to $\Q$.

To see this algebraic group concretely, fix a $\Q$-basis of~$k$.  Let
$k_1$ be any commutative $k$-algebra.  A point of $\bold G$ with
entries in $k_1$ is an $n\times n$ matrix with non-zero determinant
whose $(i,j)$-th entry is a $d$-tuple of coordinates $(x_{ij}^{(1)},
\dots, x_{ij}^{(d)})$ with $x_{ij}^{(l)} \in k_1$.  The $d$-tuple
gives an element of $k \otimes_\Q k_1$ with respect to the fixed
$\Q$-basis of~$k$.  Addition and multiplication of $d$-tuples are
defined using the corresponding operations in $k$.  We fix such a
realization of $\bold G$ throughout the paper.

\subheading{(1.3)} Let $S = k \otimes_\Q \R$. The real and complex
embeddings of $k$ induce a canonical isomorphism of $\R$-algebras $S
\cong \prod_v k_v$.  Clearly $k\hookrightarrow S$ by the canonical
diagonal embedding.

Set $G = \bold G(\R) = \GL_n(S)$.  Each entry in a matrix of $G$
breaks up as a product over the $v$'s, so that $G \cong \prod_v
\GL_n(k_v)$.

\subheading{(1.4)} Let $\qT$ be the following maximal $\Q$-split torus of
$\bold G$: if the fixed basis of $k$ over $\Q$ has~1 as its first
element, then $\qT$ is the subgroup where all the
$x_{ij}^{(m)}$ are zero except for $x_{ii}^{(1)}$, $i=1,\dots,n$.
Write $\qT = \{\op{diag}(a_1,\dots,a_n)\}$ with $a_i = x_{ii}^{(1)}$.

The split radical $R_d \bold G$ of $\bold G$ is the subgroup of $\qT$
in which all the $a_i$ are equal to each other.  Let $H$ be the
identity component of $R_d \bold G(\R)$.  These are the {\sl positive
real homotheties\/}: $h\in H$ acts by a common positive real scalar on
each factor $\GL_n(k_v)$ of $G$.  We identify $H \cong \R_+$.

\subheading{(1.5)} Let $\bold a, \bold b \in S^n$.  These break up as
a product of vectors $\bold a_v, \bold b_v \in (k_v)^n$.  Write $\bold
a_v = (a_1,\dots,a_n)$ with $a_i \in k_v$, and write $\bold b_v$
similarly.  For each $v$, let $\langle (a_1,\dots,a_n),
(b_1,\dots,b_n) \rangle_v = \sum_{i=1}^n \op{trace}^{k_v}_\R (a_i \bar
b_i)$, where the bar denotes complex conjugation when $v$ is complex
and the identity when $v$ is real.  Let $\langle \bold a, \bold b
\rangle = \sum_v \langle \bold a_v, \bold b_v \rangle_v$, a
positive-definite product on $S^n$.  Let $K_v$ be the subgroup of
$\GL_n(k_v)$ preserving $\langle \, , \, \rangle_v$.  Let $K = \prod_v
K_v$; this is the largest subgroup of $G$ preserving $\langle \, , \,
\rangle$, and is a maximal compact subgroup of $G$.  Let $\| \bold a
\| = \langle \bold a, \bold a \rangle^{1/2}$.

\definition{Definition} Let $X = (HK) \bs G$.  \enddefinition

In the language of \cite{B-S}, $X$ is a homogeneous space for $G$ of
type $S$--$\Q$.  It is a symmetric Riemannian space for the Lie
group~$H\bs G$.

\subheading{(1.6)} {\it Example of $\GL_n(\Z)$.} Let $k=\Q$.  Then
$G=\GL_n(\R)$, the group of $n\times n$ invertible matrices with real
entries.  The group $\qT$ consists of the diagonal matrices, and $H =
\{ \op{diag}(h,\dots,h) \mid h > 0 \}$.  The space $X$ is the
symmetric space consisting of the real $n\times n$ symmetric
positive-definite matrices mod homotheties.  When $n=2$, $X$ is
isomorphic to the upper half-plane (see~(2.6)).


\heading Section 2---Marked Lattices.  The Well-Rounded Retract
\endheading

We recall the definition and basic properties of the well-rounded
retract for $\GL_n$.  Except for the emphasis on marked lattices, the
material comes from~\cite{A3}; we specialize this to the number field
case, and omit the proofs.  For the $\GL_n(\Z)$ case, see~(2.6). 

\subheading{(2.1)} The space $S^n$ of column vectors is an $\R$-vector
space.  $S$ (and hence $k\subset S$) acts on it by coordinate-wise
multiplication, and $G = \GL_n(S)$ acts on it on the left.  These
module structures commute with each other.  If $A$ is a subset of
$S^n$, the $S$-module it spans in $S^n$ is denoted $S\cdot A$.

\subheading{(2.2)} A {\sl $\Z$-lattice\/} in a vector space $V$ over
$\Q$ is a finitely generated subgroup of $V$ which contains a
$\Q$-basis of $V$.  A $\Z$-lattice in a vector space over~$\R$ is a
finitely generated subgroup which contains an $\R$-basis that is also
a $\Z$-basis.

Fix once and for all a $\Z$-lattice $L_0$ in $k^n$.  Let $\cal O_0$ be
the subset of $k$ which stabilizes $L_0$.  Then $\cal O_0$ is an order
in $k$.  Let $\Gamma_0$ be the arithmetic group consisting of the
$\gamma \in G$ such that $\gamma \cdot L_0 = L_0$.  Let $\Gamma$ be an
arithmetic subgroup of finite index in~$\Gamma_0$.

We shall identify lattices that differ only by a homothety, namely $L$
and $hL$ for $h\in H$.  Let $Y$ be the set of all $\Z$-lattices in
$S^n$ that are stable under $\cal O_0$ and are isomorphic to $L_0$ as
$\cal O_0$-module.  Then we may identify $Y$ with $H\bs G/\Gamma_0$,
where for $g\in G$ the lattice $L = g L_0$ corresponds to the coset
$g\Gamma_0$.  Give $Y$ the topology coming from $H \bs G/\Gamma_0$.

For any subset $A \subseteq L = g L_0$, the $S$-span $S\cdot A$ is a
free $S$-module whose $S$-rank equals the dimension over $k$ of the
$k$-vector space spanned by $g^{-1}A$.

\subheading{(2.3)} The group $K$ acts on $Y$, and in fact
$$
K\bs Y = X/\Gamma_0,
$$
since both are $(HK)\bs G / \Gamma_0$.  However, the points of view on
the two sides of this equation are different: $K\bs Y$ is a space of
lattices modulo rotations, while $X/\Gamma_0$ is a homogeneous space
modulo an arithmetic group.  We now make explicit how to go back and
forth between the two sides.

\definition{Definition} A {\sl marked lattice in $S^n$\/} is a
function $f: L_0 \to S^n$ of the form $f(\bold x) = g \bold x$ for
$g\in G$.  We shall identify marked lattices that differ only by a
homothety, i.e.~$f$ and $h f$ for $h\in H$.  \enddefinition

Let $Y' = H\bs G$ be the set of marked lattices.  A marked lattice $f$
gives rise to an ordinary lattice $L\in Y$ by setting $L=f(L_0)$.
This realizes the projection $Y' = H\bs G \to Y = H\bs G /\Gamma_0$.
On the other hand, $f(\bold x) = g \bold x$ gives a point $(HK)g \in
X$, which realizes the projection $Y' = H\bs G \to X = (HK)\bs G$.
The right action of $g_1\in G$ on $Y'$ sends $f: \bold x \mapsto
g\bold x$ to $\bold x \mapsto g g_1 \bold x$.  The left action of $k_1
\in K$ on~$Y'$ sends $f$ to $k_1 \cdot f$.  Diagram~(2.3.1) shows
these spaces.

$$
\matrix
 &&
  \matrix Y' = H\bs G \\ \text{(marked lattices)} \endmatrix
  && \\
 &\swarrow&&\searrow& \\
 \matrix X = HK\bs G \\ \text{(homogeneous space)} \endmatrix
  &&&&
  \matrix Y = H\bs G / \Gamma_0 \\ \text{(lattices)} \endmatrix \\
 &\searrow&&\swarrow& \\
 &&
  \matrix X/\Gamma_0 = K\bs Y \\ \text{(arithmetic quotient)} \endmatrix
  &&
\endmatrix
\tag 2.3.1
$$

\subheading{(2.4)} Let $f$ be a marked lattice, fixed within its
homothety class.  Let $L = f(L_0)$.  The {\sl arithmetic minimum of\/}
$f$ is defined to be $m(f) = \min \{ \| \bold a \| \mid \bold a \in L
- \{0\} \}$.  This number is positive.  The set of {\sl minimal
vectors of\/} $f$ is defined to be $M(f) = \{ \bold a \in L \mid \|
\bold a \| = m(L) \}$.

The definitions of $m(f)$ and $M(f)$ use only the image $L$ of $f$, so
they descend to functions $m(L)$ and $M(L)$ on $Y$.  On the other
hand, $m(f)$ and $M(f)$ are $K$-equivariant, since $K$ preserves
$\|\dots\|$.

Unless otherwise specified, we normalize $f$ (resp.~$L$) in its
homothety class so that $m(f)$ (resp.~$m(L)$) equals~1.

\definition{{\bf (2.5).} Definition} A marked lattice $f$ is {\sl
well-rounded\/} if $M(f)$ spans $S^n$ as $S$-module. \enddefinition

By the $\Gamma_0$-invariance and $K$-equivariance noted in~(2.4), we
get a definition of well-rounded elements in all the spaces
in~(2.3.1).  For instance, a point in $X$ is well-rounded if and only
if, for the class of marked lattices $f$ mod~$K$ that represents it,
one (and hence each) of the $M(f)$ spans $S^n$ as $S$-module.  The
definitions in $X/\Gamma_0$ and $Y$ are similar.

\definition{Definition} $W$ is the set of well-rounded elements in
$X$.  We call this the {\sl well-rounded retract\/} in $X$.  For any
arithmetic group $\Gamma\subseteq \Gamma_0$, $W/\Gamma$ is called the
{\sl well-rounded retract\/} in $X/\Gamma$.  \enddefinition

\noindent In Section~4, we will explain why these are deformation
retracts.

\subheading{(2.6)} {\it Example of $\GL_n(\Z)$.}  As in~(1.6), take
$k=\Q$ and $G=\GL_n(\R)$.  If $L_0$ is the standard lattice $\Z^n
\subset \R^n$, then $\Gamma_0 = \GL_n(\Z)$, and $Y$ is the set of all
rank-$n$ lattices in $\R^n$ mod homotheties.  The product $\langle
\,,\, \rangle$ is the standard dot product.

Let's work out the $\GL_2(\Z)$ case in detail.  For each
marked lattice $f : L_0 \to \R^2$, rotate the image until
$f((1,0))$ lies on the positive $x$-axis, and use a homothety to
make $\| f((1,0)) \| = 1$.  This identifies $X$ with the upper
half-plane~$\frak h$ by sending $x = HK
\left(\smallmatrix 1 & a \\ 0 & b \endsmallmatrix\right)$ (with
$b>0$) to the point $a+bi \in \frak h$.  The well-rounded retract $W$ is
the tree shown in the figure below.  The lower branches break in two
infinitely many times as they approach the horizontal axis $\{b=0\}$,
which is not pictured.

\yawn{2.5truein}

The group $G$ is generated by $\SL_2(\R)$, $H$, and
$\left(\smallmatrix 1 & 0 \\ 0 & -1 \endsmallmatrix\right)$.  Under
our identification of $X$ with $\frak h$, $\SL_2(\R)$ acts on $\frak
h$ by linear fractional transformations as usual, $H$ acts trivially,
and $\left(\smallmatrix 1 & 0 \\ 0 & -1 \endsmallmatrix\right)$ acts
by $a+bi \mapsto -a+bi$.  The well-known fundamental domain for the
action of $\SL_2(\Z)$ on $\frak h$ is $F = \{a+bi \in \frak h \mid |a|
\le \frac12, a^2 + b^2 \ge 1 \}$, shown with hatching.  A fundamental
domain for $\GL_2(\Z)$ is the right-hand half of $F$, namely $F^* = \{
a+bi \in F \mid a \ge 0\}$, shown with double hatching.  The space
$X/\Gamma_0$ is homeomorphic to $F^*$.  The retract $W/\Gamma_0$ is
identified with the bottom arc of $F^*$.


\heading Section 3---The Borel-Serre Compactification \endheading

In this section we recall the construction of the Borel-Serre manifold
with corners $\bar X$, and we summarize facts about flags and
parabolic subgroups.  The main goal of the section is Proposition~3.6,
which relates the geodesic action on $X$ to the idea of rescaling a
lattice in directions given by a fixed flag~$\F$.

\subheading{(3.1)} Recall that $\qT = \{ \op{diag}(a_1,\dots,a_n) \}$
is a maximal $\Q$-split torus in $\bold G$.  The set of $\Q$-roots of
$\bold G$ with respect to $\qT$ is $\{a_i a_j^{-1} \mid 1 \le i \ne j
\le n\}$.  Fix the fundamental system of simple roots $\Delta = \{a_i
a_{i+1}^{-1} \mid i = 1,\dots,n-1 \}$.

\subheading{(3.2)} We recall some facts about parabolic $\Q$-subgroups
and flags.  Let $\{e_1,\dots, e_n\}$ be the standard $k$-basis of
$k^n$.  Let $J\subseteq \Delta$.  Form a graph with vertices in $\{1,
\dots, n\}$ and with an edge joining vertices $i$ and $i+1$ if and
only if $a_i a_{i+1}^{-1} \in J$.  The graph has $l = \#(\Delta -
J)+1$ connected components.  Say the $j$-th component in order from
left to right has $\nu_j$ vertices.  Consider the flag $\F_J = \{ 0
\subsetneqq V_1 \subsetneqq \cdots \subsetneqq V_l = k^n \}$ in which
$V_j$ is spanned by $\{e_1, \dots, e_{\nu_1+\dots+\nu_j} \}$.  The
standard parabolic $\Q$-subgroup $\bold P_J$ is the subgroup of $\bold
G$ such that $\bold P_J(\Q)$ is the stabilizer of $\cal F_J$ in $\bold
G(\Q)$.  In the coordinates of~(1.2), $\bold P_J$ has the usual block
upper-triangular form, with diagonal blocks of sizes $\nu_1\times
\nu_1, \dots, \nu_l\times\nu_l$.  Every parabolic $\Q$-subgroup $\bold
P$ of $\bold G$ is conjugate by some $g\in \bold G(\Q)$ to a unique
$\bold P_J$.

Let $P_J = \bold P_J(\R)$, and $P = \bold P(\R)$ in general. 

For $g \in \bold G(\Q)$, the group of $\Q$-points of $g^{-1} \bold P_J
g$ is the stabilizer of the flag $\cal F = g^{-1} \cal F_J$.  Any flag
of the latter form is called a {\sl $\Q$-flag}.  There is a one-to-one
correspondence between the parabolic $\Q$-subgroups and the $\Q$-flags
they stabilize.  We write $\F' \supseteq \F$ if every $V_j$ belonging
to $\F$ is also a member of $\F'$.

For the $\Q$-flag $\F = \{ 0 \subsetneqq V_1 \subsetneqq \cdots
\subsetneqq V_l = k^n\}$, we set $\#\F = l$.  Thus $\#\F = n$ when $P$
is a Borel subgroup, $\#\F = 2$ for maximal proper $P$, and $\#\F = 1$
when $P = G$.

\subheading{(3.3)} Let $\qT_J$ be the subgroup of $\qT$ in which the
$i$-th and $(i+1)$-st diagonal entries are equal whenever $t_i
t_{i+1}^{-1} \in J$.  Let $T_J$ be the identity component of
$\qT_J(\R)$, and let $A_J = H\bs T_J$.

Let $P = g^{-1} P_J g$ for $g\in \bold G(\Q)$.  There is a semi-direct
product decomposition $P = M_P T_P U_P$.  Here $U_P$ is the
uni\-potent radical of $P$, $T_P = g^{-1} T_J g$, and $M_P T_P$ is the
unique Levi subgroup of $P$ which is stable under the Cartan
involution of $G$ which fixes $g^{-1} K g$.  We set $A_P = H\bs T_P$.

When $P \subseteq P'$, $A_{P'}$ is naturally identified
with a subgroup of $A_P$ by taking the kernel of appropriate roots
in~$\Delta$.

\subheading{(3.4)} We recall the definition of the geodesic action
\cite{B-S, \S3}.  Let $P = g^{-1} P_J g$ for $g\in \bold
G(\Q)$.  Let $Z$ be the identity component of the center of $M_P T_P$.
Let $y\in X$.  There is a $p\in P$ with $y = HKgp$.  For any $z \in Z$,
define
$$
y \geod z = HKgzp.
$$
The right side is independent of the choices of $g$ and $p$.  This
action of $Z$ is the {\sl geodesic action\/} on~$X$ associated
to $P$.  It is a proper and free action, and it commutes
with the ordinary action of~$P$.  The group $A_P$ is contained in
$H\bs Z$, so it operates on $X$ by the geodesic action.

\subheading{(3.5)} We recall the construction of $\bar X$ \cite{B-S,
\S\S5--7}.  The roots in $\Delta - J$ determine an isomorphism between
$A_J$ and $(0,\infty)^{l-1}$, sending the point
$$
\op{diag}(\undersetbrace{\nu_1}\to{a_1,\dots,a_1}, \dots,
\undersetbrace{\nu_l}\to{a_l,\dots,a_l})
\qquad (a_j > 0, \text{ mod }H)
$$
to $(a_1 a_2^{-1}, \dots, a_{l-1} a_l^{-1}) \in (0,\infty)^{l-1}$.
Let $\bar A_J$ be the partial compactification of $A_J$ corresponding
to $[0,\infty)^{l-1}$ under this isomorphism.  The {\sl corner\/}
$X(P_J)$ {\sl associated to $P_J$} is the fiber product
$X(P_J) = X \times^{A_J} \bar A_J$, where $A_J$ acts on $X$ by
the geodesic action.  Conjugating by an element of $\bold G(\Q)$, we
define $\bar A_P$ and a corner $X(P) = X \times^{A_P} \bar A_P$
for any~$P$.

Whenever $P \subseteq P'$, there is a map $X(P') \hookrightarrow X(P)$
induced by the inclusion on $X$; we identify $X(P')$ with its image
under this map.  As a set, $\bar X$ is the union of the $X(P)$ over
all $P$.  One puts on $\bar X$ the topology such that the original
topology on each corner is preserved, $\bar X$ is Hausdorff, and
$\Gamma$ acts properly on $\bar X$ with a compact quotient.

Let $e(P)$ be the set corresponding to $X \times^{A_P} \{0\}^{l-1}$
in the corner $X(P) = X \times^{A_P} \bar A_P$.  As a set, $\bar
X$ is the disjoint union of $X$ and the $e(P)$ for all proper
parabolic $\Q$-subgroups $P$.  We call $e(P)$ the {\sl
boundary face corresponding to\/} $P$.

Let $\bar e(P)$ be the closure of $e(P)$ in $\bar X$.  We have $\bar
e(P) = \coprod_{P' \subseteq P} e(P')$.  The sets $\bar e(P_1)$ and
$\bar e(P_2)$ are disjoint unless $P_1 \cap P_2$ is a parabolic
$\Q$-subgroup~$P_3$, in which case $\bar e(P_1) \cap \bar e(P_2) =
\bar e(P_3)$.

Notice that the largest subgroup of $\Gamma$ which stabilizes $X(P)$,
$e(P)$, and $\bar e(P)$ is $\Gamma\cap P$.

Modding out by the geodesic action of $A_P$ gives a $(\Gamma\cap
P)$-equivariant bundle map $q_P : X(P) \to e(P)$.  The
orbits in $X(P)$ under $M_P U_P$ are sections of this bundle
called the {\sl canonical cross-sections.}

\subheading{(3.6)} We now present a central tool of the paper, an
interpretation of the geodesic action in the setting of marked
lattices.

Let $f$ be a marked lattice.  Let $\cal F = \{0 = V_0 \subsetneqq
\cdots \subsetneqq V_l = k^n \}$ be a $\Q$-flag, with $P = g^{-1} P_J
g$ the parabolic $\Q$-subgroup determined by~$\F$.  Throughout the
paper, we will abuse notation by writing $f(V_j \cap L_0) \otimes_{\Q}
\R \subseteq S^n$ as $f(V_j)$; this is an $S$-linear subspace of
$S^n$, of $S$-rank equal to $\dim_k V_j$.  Set $\check V_j$ equal to
the orthogonal complement of $f(V_{j-1})$ in $f(V_j)$ with respect to 
$\langle,\rangle$.  We have $S^n = \check V_1 \oplus \cdots \oplus
\check V_l$ as an orthogonal direct sum.

\definition{Definition} Let $\frak A_{\F,f}^*$ be the group of
$\R$-linear maps $S^n \to S^n$ which act by positive real homotheties
on each summand $\check V_j$ separately.  Let $\frak A_{\F,f} = \frak
A_{\F,f}^* / H$.  We call $\frak A_{\F,f}$ the {\sl orthogonal scaling
group with respect to~$\F$ and~$f$}. \enddefinition

Note that every $\alpha \in \frak A_{\F,f}^*$ is $S$-linear.  The
group $\frak A_{\F,f}$ acts on $Y'$ via $f_1 \mapsto \alpha \cdot
f_1$.  However, it does not act on~$X$.  (The next proposition will
show that $\frak A_{\F,f}$ acts on $Y' = H\bs G$ as a lift of the
action of $A_P$ on~$X = HK\bs G$, but this lift is only determined by
an arbitrary choice of an element of~$K$.  The choice is encoded
in~$f$.)

Let $\Psi : Y' = H \bs G \to X = HK \bs G$ be the projection.  Let $y
= \Psi(f)$.  Let $\iota : A_P \to \frak A_{\F,f}$ be the following
map: if $z \in A_P$ with
$$
z = g^{-1} \op{diag}(\undersetbrace{\nu_1}\to{a_1,\dots,a_1}, \dots,
\undersetbrace{\nu_l}\to{a_l,\dots,a_l}) g  \qquad (a_j > 0), \tag 3.6.1
$$
let $\iota(z)$ be the element of $\frak A_{\F,f}$ which (for each
$j=1,\dots,l$) acts on $\check V_j$ by the scalar~$a_j$.  It is clear
that $\iota$ is an isomorphism of groups.

\proclaim{Proposition} For any $z\in A_P$, with $\iota$ as above,
$$
y \geod z = \Psi(\iota(z) \cdot f). \tag 3.6.2
$$
\endproclaim

\remark{Remark} The proposition is stated in \cite{G, (2.4)}, without
proof.  \endremark

\demo{Proof} We use the notation of~(3.4).  Write $y = HKgp$ and $p =
g^{-1} p_J g$ for $p_J \in P_J$.  Consider the marked lattice $f_1 :
\bold x \mapsto gp\bold x = p_J g \bold x$.  Since $\Psi(f_1) = y =
\Psi(f)$, we may choose a $k_1 \in K$ such that $f_1 = k_1 \cdot f$.

Write $z\in A_P$ as $z = g^{-1} z_J g$ for
$z_J = \op{diag}(a_1,\dots,a_1,\dots, a_l,\dots,a_l) \in A_J$.  Then
$$
y\geod z = HK gzp = HK z_J p_J g.
$$
The marked lattice $f_z : \bold x \to z_J p_J g \bold x$ descends
mod~$K$ to $y\geod z$.

Let $\cal I_J$ be the flag in $S^n$ whose $j$-th member is the
coordinate subspace spanned by $\{ e_1, \dots, e_{\nu_1 + \cdots +
\nu_j} \}$.  (This $\cal I_J$ is essentially $\F_J \otimes_{\Q} \R$.)
The map $\bold x \mapsto g\bold x$ carries $\F$ to $\cal I_J$, and the
further multiplication by $p_J$ and $z_J$ preserves $\cal I_J$.  Thus
for any $z\in A_P$, $f_z(V_j)$ is exactly the $j$-th member of $\cal
I_J$.  (The idea is that, out of the whole $K$-equivalence class of
marked lattices representing $y\geod z$, $f_z$ is a representative for
which the $f_z(V_j)$ will be in standard position.)  The
orthocomplement of $f_z(V_{j-1})$ in $f_z(V_j)$ is the coordinate
subspace $E_j$ spanned by $\{e_{\nu_1 + \cdots + \nu_{j-1} + 1}, \dots,
e_{\nu_1 + \cdots + \nu_j}\}$.  The $E_j$ are orthogonal, and $z_J$
acts by the map $\alpha_1 : S^n \to S^n$ characterized by the property
that it acts on each $E_j$ by the scalar $a_j$.  In short, the
left-hand side of~(3.6.2) is the projection under $\Psi$ of the marked
lattice given by the composition $L_0 \overset{f_1} \to \to S^n
\overset{\alpha_1} \to \to S^n$.

On the right-hand side of~(3.6.2), observe that $\Psi(\alpha\cdot f) =
\Psi(k_1 \alpha \cdot f) = \Psi((k_1 \alpha k_1^{-1}) f_1)$.  Since
$f_1 = k_1 \cdot f$, the orthogonal complement of $k_1 \cdot
f(V_{j-1})$ in $k_1 \cdot f(V)$ equals that of $f_1(V_{j-1})$ in
$f_1(V_j)$; this says exactly that $k_1 \cdot \check V_j = E_j$ for
all~$j$.  It is immediate that $k_1 \alpha k_1^{-1}$ is the same map
as~$\alpha_1$.  Thus the right-hand side of~(3.6.2) is also $\Psi$
applied to the composition $L_0 \overset{f_1} \to \to S^n
\overset{\alpha_1} \to \to S^n$. \qed\enddemo

\remark{Note} We will always normalize elements of $A_P$ and $\frak
A_{\F,f}$ in their classes modulo~$H$ so that $a_1 = 1$.
\endremark

\subheading{(3.7)} While $A_P$ and $\frak A_{\F,f}$ act on different
spaces, we use similar notation for the two actions.  We write 
$$
\pmb\rho = (\rho_1, \dots, \rho_{l-1}) \in (0,\infty)^{l-1}
$$
for a point in either group.  The~$\rho_i$ are the
coordinates  coming from the simple roots in~$J$: that
is, $\rho_j = a_j /
a_{j+1}$ for $j=1,\dots,l-1$.
By Proposition~3.6, the point of $\frak A_{\F,f}$ corresponding to $(\rho_1,
\dots, \rho_{l-1})$ acts on $\check V_j$ by the scalar
$$
a_j = (\rho_1 \cdots \rho_{j-1})^{-1}  \qquad(j=2,\dots,l). \tag 3.7.1
$$
(Because of our normalizations, the action on
$\check V_1$ is by~1.)

\remark{Remark} The $\rho_j$'s are good coordinates for rescaling
sublattices within lattices, while the $a_j$ are the natural
coordinates for the geodesic action.  Though~(3.7.1) is messy, we
cannot avoid using it, since part of our goal is to relate lattices
and the geodesic action.  The inverses in~(3.7.1) arise because the
well-rounded retraction goes away from the Borel-Serre
boundary. \endremark

\smallskip

Multiplication in $A_P$ corresponds to coordinate-wise multiplication
of the $\rho_i$.  As in~(3.5), the $\rho_j$ extend to coordinates on
$\bar A_P$ with $(\rho_1, \dots, \rho_{l-1}) \in
[0,\infty)^{l-1}$.  Going to infinity in $X(P)$
means going to zero in the $\rho_j$ variables.  We introduce a partial ordering on $A_P$ and $\bar
A_P$, where $(\rho_1, \dots, \rho_{l-1}) \le (\rho_1', \dots,
\rho_{l-1}')$ if and only if $\rho_j \le \rho_j'$ for each~$j$.

The action of $\pmb\rho = ( \rho_1, \dots, \rho_{l-1})$ on $Y'$ is
denoted $f \mapsto \pmb\rho \cdot f$.  This descends mod~$\Gamma_0$ to
an action $L \mapsto \pmb\rho \cdot L$ on~$Y$.  Descending mod~$K$ as
in Proposition~3.6, $\pmb\rho \cdot x$ denotes the geodesic action on
points $x$ of $X$ or $X/(\Gamma\cap P)$.

\subheading{(3.8)} We will need to consider $\Q$-flags $\cal M = \{0
\subsetneqq \Ms 1 \subseteq \Ms 2 \subseteq \cdots \subseteq \Ms n =
k^n \}$ in which some of the inclusions can be equalities.  If $\F$ is
obtained from $\cal M$ by replacing each string of equalities with a
single space $V_j$, we call $\F$ the {\sl irredundant version\/} of
$\cal M$.  If $\frak A$ is the orthogonal scaling group with respect
to this flag and some~$f$, we always write $\pmb\rho = (\rho_1, \dots,
\rho_{l-1}) \in \frak A$ with reference to $\F$, not $\cal M$.


\heading Section 4---The Well-Rounded Retraction \endheading

In this section, we define the well-rounded retraction $r_t$.  The
material comes from~\cite{A3}, though we emphasize different things
(marked lattices, the flag of successive minima).

\subheading{(4.1)} Let $f$ be a marked lattice, with $L = f(L_0)$.  We
normalize $f$ up to homothety so that $m(f) = 1$.  If $V\subseteq k^n$
is a $k$-subspace, $f(V) \otimes_{\Q} \R$ is abbreviated $f(V)$.

Recall that $Y' = H\bs G$ is the space of marked lattices~(2.3).  For
$i=1,\dots,n$, let
$$
Y'_i = \{ f\in Y' \mid \op{rank}_S (S\cdot M(f)) \ge i\}.
$$
Observe that $Y' = Y'_1 \supset Y'_2 \supset \cdots$, and
that $Y'_n$ is the set of well-rounded elements of $Y'$. 
For $i=1,\dots,n-1$, we will define a
deformation retraction ${\rts i} : Y'_i \times [0,1] \to Y'_i$ with
image $Y'_{i+1}$.  The deformation
retraction will be equivariant for both the $K$ and $\Gamma_0$
actions.  For any $f\in Y'$, we will also define a
$\Q$-flag $\cal M = \{ 0 \subseteq \Ms 1 \subseteq \cdots \subseteq
\Ms n = k^n \}$ depending on $f$.  (This $\cal M$ could have
equalities, as in~(3.8).)  The definitions proceed simultaneously by
induction on~$i$.

Let $f_1 = f$, and for $i>1$ let $f_i = r^{(i-1)}_1(f_{i-1})$ inductively.
Let $\Ms i$ be the $k$-span of $f_i^{-1}(M(f_i))$.  This is a subspace of
$k^n$ of dimension~$\ge i$.  If $f_i$ is in $Y'_{i+1}$ already, set
${\rts i}(f_i) = f_i$ for all $t\in [0,1]$, and set $\mu_i = 1$.
Otherwise, consider the $S$-linear maps $\phi_\mu$ obtained by
multiplying in $f_i(\Ms i)$ by~1 and in the orthocomplement $f_i(\Ms
i)^\perp$ by a scalar $\mu>0$.  
By \cite{A3, p.~463}, there is a unique smallest $\mu_i \in (0,1)$
such that
$$
\phi_{\mu_i} \cdot f_i
 \text{ still has arithmetic minimum equal to 1.}
\tag 4.1.1
$$
This $\phi_{\mu_i} \cdot f_i$ lies in $Y'_{i+1}$.  We set ${\rts
i}(f_i) = \phi_{1+(\mu_i-1)t} \cdot f_i$.
This is continuous in both the $f_i$ and $t$ variables, and defines a
homotopy between $r_0^{(i)} = (\text{identity on $Y_i'$})$ and
$r_1^{(i)} : Y'_i \twoheadrightarrow Y'_{i+1}$.

\definition{Definition} The {\sl well-rounded retraction\/} $r_t$ on
$Y'$ is
obtained by composing the deformation retractions $r_t^{(i)}$ in the
order $i=1,\dots,n-1$.  More precisely, $r_t : Y'\times [0,1] \to Y'$ is
given by
$$
r_t(f) = r^{(i)}_{t(n-1)-(i-1)} \circ r_1^{(i-1)} \circ \cdots \circ
                        r_1^{(1)}(f),
$$
where $i$ is determined by the rule $t \in [ \frac{i-1}{n-1},
\frac{i}{n-1} ]$ for $i \in \{1,\dots,n-1\}$.

Let $r = r_1$.  We will also call this the {\sl well-rounded
retraction}. \enddefinition

Notice that $\Ms {i+1} \supseteq \Ms i$ by construction.  This gives the

\definition{Definition} The {\sl flag of successive minima\/} for $f$
is the $\Q$-flag $\cal M = \{ 0 \subseteq \Ms 1 \subseteq \dots
\subseteq \dots \Ms n = k^n \}$.  \enddefinition

This $\cal M$ is an invariant of $f$'s class mod~$K$, and is
equivariant for $\Gamma_0$.  It is a ``lower semi-continuous''
function of $f$ in the sense that for any $f \in Y'$, there is a
neighborhood $U$ of $f$ such that for all $f' \in U$, the flag $\cal
M'$ of successive minima for $f'$ satisfies $\cal M' \subseteq \cal
M$.

\remark{{\bf (4.2).} Remark} Once the inductive construction of $r$ is
completed, we obtain the following interpretation of $r(f)$ as the
result of applying a geodesic action to~$f$.  Let $\cal M$ be the flag
of successive minima for $f$, with irredundant version~$\F$.  Let
$\pmb\mu = (\mu_1^{-1},
\dots, \mu_{n-1}^{-1}) \in \frak A_{\F,f}$ with the $\mu_i$ depending on~$f$
as in~(4.1).  (As in~(3.8), $\mu_i^{-1}$ is omitted when $\Ms i = \Ms
{i+1}$.)  Then
$$
r(f) = \pmb\mu \cdot f.
$$
\endremark

\subheading{(4.3)} We now adapt the definition of the well-rounded
retraction to the other spaces in~(2.3.1).  In each case, the
retraction will still be denoted $r_t$ for $t\in [0,1]$, with $r=r_1$,
and will be called the ``well-rounded retraction''.

The whole construction is $\Gamma_0$- and $K$-equivariant.  Hence
$r_t$ descends to deformation retractions of $Y$ onto the set of
well-rounded lattices, of $X$ onto $W$, and of $X/\Gamma$ onto
$W/\Gamma$ for any $\Gamma\subseteq \Gamma_0$ of finite index.

We recall the main theorem of \cite{A3}.

\proclaim{Theorem} The map $r_t$ is a $\Gamma_0$-equivariant
deformation retraction of $X$ onto $W$.  The quotient $W/\Gamma_0$ is
compact.  The dimension of $W$ and $W/\Gamma_0$ is the virtual
cohomological dimension of~$\Gamma_0$, which equals $\dim X - (n-1)$.
The space $W$ has a natural structure as a cell complex on which
$\Gamma_0$ acts cell-wise with finite stabilizers of cells.  The first
barycentric subdivision of the cell structure descends to a finite
cell complex structure on $W/\Gamma_0$. \endproclaim

The example of $\GL_2(\Z)$ is discussed in~(6.3).

\remark{Remark} When the class number $h$ of $k$ is~$>1$, \cite{A3}
defines a whole $(h-1)$-dimensional family of retracts, depending on a
parameter called a {\sl set of weights}.  For simplicity, we will use
only the trivial set of weights (identically equal to~1).  Readers may
easily make the needed changes if they wish to use a non-trivial set
of weights. \endremark


\heading Section 5---Summary of Notational Conventions \endheading

The letter $f$ always stands for a marked lattice $f : \bold x \mapsto
g \bold x$ for $g\in G$, modulo the homotheties~$H$.  We let $L =
f(L_0)$, sometimes without explicit mention.  We let $L'$ be the
lattice for $f'$, etc.  As in~(3.6), $f(V_j\cap L_0) \otimes_\Q \R
\subseteq S^n$ will be written $f(V_j)$.  Let $\check V_j$ be the
orthogonal complement of $f(V_{j-1})$ in $f(V_j)$, so that $S^n =
\check V_1 \oplus \cdots \oplus \check V_l$ as an orthogonal direct
sum.  As in~(2.4), lattices are always scaled by a homothety so that
their shortest non-zero vector has length~1.

For the rest of the paper, $P$ denotes the real points of a parabolic
$\Q$-subgroup.  Let $\F = \{ 0 \subsetneqq V_1 \subsetneqq \dots
\subsetneqq \dots V_l = k^n \}$ be the $\Q$-flag for $P$ as in~(3.2),
and $\frak A = \frak A_{\F,f}$ the orthogonal scaling group with
respect to~$\F$ and~$f$ (see (3.6)--(3.7)).

{\sl Whenever we mention $P$, the letters $\F$ and $\frak A$
will refer to the objects corresponding to $P$, and vice versa.}
The symbol $P'$ is associated to $\F'$ and $\frak A'$, etc.

{\sl From Section~6 on, we will often write $f$ for the point that $f$
determines in $X$ or $X/\Gamma$}.  We will speak of marked lattices
$f$ ``in'' $X$ or $W$, or write $L$ for the point that $L$ determines
in $X/\Gamma$.  These abuses are justified because our constructions
involving marked lattices are $K$-equivariant.


\heading Section 6---The Sets $W_\F$ \endheading

\subheading{(6.1)} Throughout this section, $\cal F = \{0 = V_0
\subsetneqq V_1 \subsetneqq \cdots V_l = k^n \}$ will be a fixed
$\Q$-flag.

\definition{Definition} Let $W_{\cal F}$ be the set of all
marked lattices $f$ in $W$ such that for all $j=1,\dots,l$,
$$
S\cdot (M(L) \cap f(V_j)) = f(V_j).
$$
\enddefinition

($W_\F$ was called $W'$ in the introduction.) Geometrically, $W_\F$
consists of the lattices $L=f(L_0)$ whose intersection with each
subspace $f(V_j)$ is well-rounded in that subspace.  By the definition
of the cell structure on $W$ \cite{A3, pp.~466--7}, $W_{\cal F}$ is a
closed subcomplex of~$W$ and is stable under $\Gamma \cap P$.  A cell
in $W_{\cal F}$ is said to {\sl respect $\cal F$}.  If $\cal F =
\gamma \cal F'$ for $\gamma \in \Gamma$, then the images of $W_{\cal
F}$ and $W_{\cal F'}$ in $W/\Gamma$ are equal.  The $\Q$-flags fall
into only finitely many equivalence classes mod~$\Gamma$, so only
finitely many distinct subcomplexes of $W/\Gamma$ arise in this way.

If $\cal M$ is a $\Q$-flag with possible equalities~(3.8), and $\F$ is
the irredundant version of $\cal M$, let $W_{\cal M} = W_{\F}$.

\subheading{(6.2)}  The following is immediate from the
definitions in~(4.1).

\proclaim{Lemma} If $f$ has flag of successive minima $\cal M$, then
$r(f) \in W_{\cal M}$. \endproclaim

\subheading{(6.3)} {\sl Example of $\GL_2(\Z)$.} Refer again to the
picture in~(2.6).  For the standard flag $\cal F = \{ 0 \subsetneqq
\Q\cdot(1,0) \subsetneqq \Q^2\}$ with associated $P$, the space
$W_{\cal F}$ is the horizontal sequence of arcs forming the top
of~$W$.  The boundary face $e(P)$ in $\bar X$ (not shown) is the
horizontal line at $b=\infty$.  The geodesic action by $A_P$ flows
straight up, pushing $X$ upward to converge to $e(P)$.  The
well-rounded retraction pulls the region above $W$ straight down onto
$W_{\cal F}$.  For torsion-free $\Gamma$, $e(P)/(\Gamma\cap P)$ and
$W_\F/(\Gamma\cap P)$ are homeomorphic to circles, and the
well-rounded retraction induces the obvious isomorphism between the
cohomology of these spaces.  (The situation is more complicated for
$n>2$, when the boundary faces are no longer disjoint.)

\remark{Remark} For small~$n$, the $W_{\cal F} / (\Gamma \cap P)$ are
all readily computable.  See~(10.6). \endremark


\heading Section 7---Neighborhoods of Infinity in $\bar X$ \endheading

We now prove our main results explaining the connection between the
well-rounded retraction $r$, the geodesic action, and the topology of
the Borel-Serre boundary.  Parts (7.1)--(7.4) give a chain of
propositions that allow us to define a tubular neighborhood $\bar
N_{\cal F}$ of $e(P)$ (called $\bar N$ and $e$ in the introduction).
In (7.5)--(7.6) we prove that $r:X\to W$ has the continuous extension
$\bar r:\bar X\to W$, where (for all~$\F$) $\bar r$ is constant on the
fibers of $\bar N_\F$ and carries $\bar N_\F$ onto $W_\F$.

\subheading{(7.1)} Fix a marked lattice $f$, with $L = f(L_0)$.  As
in~(4.1), let $\cal M = \{0 \subseteq \Ms 1 \subseteq \dots \subseteq
\Ms n = k^n \}$ be the flag of successive minima for $f$.  Let $\F$ be
the irredundant version of $\cal M$, with associated $P$.  Let
$\frak A = \frak A_{\F,f}$.  Let $f'' = r(f)$ be the image of $f$ under the
well-rounded retraction.

The idea of the following lemma is that if $f$ is close enough to
$e(P)$, then the well-rounded retraction carries the whole orthant
$\{\pmb\rho \in \frak A \mid \pmb\rho \le 1\} \cdot f$ to a single
point of $W$.  The lemma says that, to know that $f$ is ``close
enough'' to $e(P)$, it suffices to know that $f$'s flag of successive
minima is the flag~$\F$ corresponding to $P$ (or at least contains
this~$\F$---see~(7.2)).

\proclaim{Lemma} Let $f$, $\cal M$, $\F$, $P$, and $\frak A$ be as
above.  For any $\pmb\rho = (\rho_1, \dots, \rho_{l-1}) \in \frak A$
with $\pmb\rho \le (1,\dots,1)$, $f' = \pmb\rho \cdot f$
satisfies $r(f') = r(f) = f''$.  The common image $f''$ will lie in $W_{\cal
M} = W_{\F}$.  \endproclaim

\demo{Proof} It suffices to prove the lemma when $\pmb\rho =
(1,\dots, 1, \rho_j, 1, \dots, 1)$ for a fixed $\rho_j \le 1$; the
general case follows by applying this for each~$j$.  Let~$i$ be the
unique index such that $V_j = \Ms i$ and $\Ms i \subsetneqq \Ms
{i+1}$.

The marked lattices $f$ and $f'$ agree on $L_0 \cap \Ms i$, since
$\pmb\rho$ acts by~1 in
that region.  For $\bold x \in L_0$ with $\bold x \notin \Ms i$, no vector
$f'(\bold x)$ is shorter than the corresponding vector $f(\bold x)$, since
$\pmb\rho$ magnifies lengths in the directions perpendicular to
$f(\Ms i)$ by $\rho_j^{-1} \ge 1$.  And $f(\bold x)$ cannot be a
minimal vector of $f$, since the minimal vectors lie
in $f(\Ms 1)$, which is contained in $f(\Ms i)$.  It follows that the
well-rounded 
retractions for $f$ and $f'$ will proceed exactly the same up through
the end of the $\rs {i-1}$ step---that is, if $\mu_\iota$ and $\mu_\iota'$
(for $f$ and $f'$ resp.) are as in~(4.1), then
$$
\mu_\iota = \mu_\iota' \text{ for all $\iota<i$}.
\tag 7.1.1
$$
Set $g = r^{(i-1)}_1 \circ \cdots \circ r^{(1)}_1 (f)$ and $g' =
r^{(i-1)}_1 \circ \cdots \circ r^{(1)}_1 (f')$.  It follows
from~(3.7.1) and~(7.1.1) that $g' = (1, \dots, 1, \rho_j^{-1}, 1,
\dots, 1) \cdot
g$.  The $\rts i$ stage of the retraction will, by definition, carry $g$ to
$$
(1, \dots, 1, \mu_i^{-1}, 1, \dots, 1) \cdot g,
\tag 7.1.2
$$
and carry $g'$ to
$$
(1, \dots, 1, {\mu_i'}^{-1}, 1, \dots, 1) \cdot g'.
\tag 7.1.3
$$
However, by the uniqueness statement in~(4.1.1), $\mu_i^{-1}$ must
equal ${\mu_i'}^{-1} \rho_i$.  Thus the marked lattices (7.1.2) and
(7.1.3) will coincide.  The $\rs {i+1}$ and later stages of the
retraction will agree on this common lattice, so the final values
$r(f)$ and $r(f')$ are equal.  By~(6.2), $f'' \in W_{\cal M}$.  \qed\enddemo

\subheading{(7.2)} We record a corollary of the previous lemma.  The
idea is that if $f$'s flag of successive minima refines the flag $\F'$
for some $P'$, then $f$ is ``close enough'' to $e(P')$.  After all, if
$f$ is close enough to $e(P)$, it should be close enough to any larger
Borel-Serre boundary component $e(P')$ that has $e(P)$ in its closure.

\proclaim{Corollary} Let the notation be as in~(7.1).  Let $\F'$ be any
$\Q$-flag with $\F' \subseteq \F$, and let $\frak A'$ be the
orthogonal scaling group for $\F'$ and~$f$.  Then for any $\pmb\rho \in
\frak A'$ with $\pmb\rho \le (1,\dots,1)$, $\pmb\rho \cdot f$ and $f$
map to the same point under the well-rounded retraction.  This point
lies in $W_{\F'}$. \endproclaim

\demo{Proof} This is trivial, since $\frak A'$ is naturally a
subgroup of $\frak A$, and $W_\F \subseteq W_{\F'}$. \qed\enddemo

\subheading{(7.3)} In~(7.1) we took an $f$ close to $e(P)$ and studied
the orthant $\{\pmb\rho \le 1\} \cdot f$.  Now, instead, we take an
arbitrary $f$.  We show in~(7.3) (actually in the corollary) that
there is a $\bold t\in \frak A$ which pushes $f$ so close to $e(P)$
that the orthant $\{\pmb\rho \le \bold t\} \cdot f$ still has the
property of~(7.1)---namely, $r$ maps it to a single point in $W_\F$.

\proclaim{Lemma} Choose any marked lattice $f$.  Fix a $\Q$-flag $\F =
\{ 0 \subsetneqq V_1 \subsetneqq \cdots \subsetneqq V_l = k^n \}$,
and let  $\frak A = \frak A_{\F,f}$.  Then there exists $\bold t =
(t_1,\dots, t_{l-1}) \in \frak A$, depending on $f$, such that for any
$\pmb\rho \in \frak A$ with $\pmb\rho \le \bold t$, the flag $\cal
M'$ of successive minima for $f' = \pmb\rho \cdot f$ satisfies $\cal
M' \supseteq \F$. \endproclaim

\demo{Proof} Let $j\in \{1,\dots,l-1\}$.  Let $d_j = \dim_k V_j$.  Let
$\overset\,\circ\to L = f(V_j \cap L_0)$.  This is a $\Z$-lattice in
the $S$-submodule $\overset\,\circ\to V = f(V_j \otimes_\Q \R)
\subseteq S^n$.  Perform the well-rounded retraction on
$\overset\,\circ\to L$ within $\overset\,\circ\to V$.  As in~(4.1),
this involves multiplying by constants $\mu_1, \dots, \mu_{d_j-1} \le
1$ in various subspaces of $\overset\,\circ\to V$.  (We must first
rescale $f$ by a homothety so that $m(\overset\,\circ\to L) = 1$; this
does not affect the value of the $\mu_i$.) Let $\beta_j = \mu_1 \cdot
\dots \cdot \mu_{d_j-1}$.

The orthogonal projection of $L = f(L_0)$ onto $f(V_j)^\perp$ is a lattice
$L^\dag$ under the restriction of the inner product to $f(V_j)^\perp$.
Let $\alpha_j > 0 $ be the length of the shortest non-zero vector of
$L^\dag$.

Let $t_1 = \min(1, \frac12 \alpha_1 \beta_1)$, and inductively let
$t_j = \min(1, \frac12 (t_1 \cdot \dots \cdot t_{j-1})^{-1} \alpha_j
\beta_j)$.  This means $(t_1\cdots t_j)^{-1} \ge 2 \alpha_j^{-1}
\beta_j^{-1}$ for $j=1,\dots, l-1$.

Now let $f' = \pmb\rho \cdot f$ for $\pmb\rho \le \bold t$.  Fix $j$.
Let $\bold a$ be any vector of $L$ lying outside $f(V_j \cap L_0)$.
Let $\bold b = \pmb\rho\cdot \bold a$ be the corresponding vector for
$f'$.  Notice that $\pmb\rho$ multiplies $\check V_{j+1}$ by $(\rho_1
\cdots \rho_j)^{-1}$, by~(3.7.1).  Since $\rho_j \le 1$ for
all $j$, $\pmb\rho$ multiplies by factors at least as large as
$(\rho_1 \cdots \rho_j)^{-1}$ in each summand of $\check V_{j+1}
\oplus \cdots \oplus \check V_l = f(V_j)^\perp$.

Imagine performing the well-rounded retraction on $\overset\,\circ\to
L$, but carrying out all the retractions on $L' = f'(L_0)$ itself;
more precisely, consider the lattice
$$
(1, \dots, 1, \mu_{d_j-1}^{-1}, 1, \dots, 1)
 \cdot (1, \dots, 1, \mu_{d_j-2}^{-1}, 1, \dots, 1)
 \cdot \dots \cdot (\mu_1^{-1}, 1, \dots, 1) \cdot L'
\tag 7.3.1
$$
where the $\mu_i$ are the constants for $\overset\,\circ\to L$.
(The $\mu_i^{-1}$ is in the $\iota$-th position
when the minimal vectors of $\overset\,\circ\to L$ at the $i$-th step of the retraction
span the $\iota$-th member of the irredundant version of the flag of
successive minima 
for $\overset\,\circ\to L$.)  Let $\bbs j$ be the vector corresponding to $\bold
b$ in~(7.3.1).  Observe, using~(3.7.1), that
$$
\align
\| \bbs j \|
 &\ge \| \text{(orthogonal projection of $\bbs j$ onto $f(V_j)^\perp$)} \| \\
 & = \mu_1 \cdot \dots \cdot \mu_{d_j-1} \cdot
        \| \text{(orthogonal projection of $\bold
                    b$ onto $f(V_j)^\perp$)} \|   \\
 & = \beta_j \cdot
        \| \text{(orthogonal projection of $\bold
                    b$ onto $f(V_j)^\perp$)} \|   \\
 & \ge \beta_j \cdot (\rho_1 \cdot \dots \cdot \rho_j)^{-1}
       \cdot \| \text{(orthogonal projection of $\bold
                    a$ onto $f(V_j)^\perp$)} \|   \\
 & \ge \beta_j \cdot (t_1 \cdot \dots \cdot t_j)^{-1} \cdot \alpha_j \\
 &\ge 2 \\
 &> 1.
\endalign
$$
This says that the image of $L'$ under the well-rounded retraction
on~$L$ up through the $\rs {d_j-1}$ stage is equal to~(7.3.1): no
vector outside $\overset\,\circ\to L$ will affect the stopping times of the
retraction on~$L'$, since each will have length~$> 1$ during the entire
process leading to~(7.3.1).  In particular, $f(V_j)$ is a member of
$\cal M'$, the flag of successive minima of $L'$, which is what we
wanted to prove.  \qed\enddemo

\proclaim{Corollary} With notation as in the lemma, the points
$\pmb\rho \cdot f$ with $\pmb\rho \le \bold t$ map to a common point
of $W_\F$
under the well-rounded retraction~$r$.  \endproclaim

\demo{Proof} Apply Corollary~7.2 to the marked lattice $\bold t \cdot
f$ of the lemma.  \qed\enddemo

\subheading{(7.4)} We can now construct neighborhoods of the
Borel-Serre boundary component $e(P)$ that are well behaved with
respect to the well-rounded retraction.  Fix a $\Q$-flag $\F$, with
its associated $P$.

Recall that the corner $X(P)$ is a fiber bundle with base
$e(P)$ and fiber $\bar A_P$.  By a {\sl neighborhood in
$X$ near $e(P)$\/}, we mean an open set in $X$ whose
intersection with every fiber contains $\{a \in A_P \mid a < a_0\}$
for some~$a_0 \in A_P$ depending on the fiber.

\proclaim{Proposition} There is a neighborhood $N_\F$ in $X$ near
$e(P)$ such that, on the intersection of $N_\F$ with a
fiber of $X(P)$, the well-rounded retraction $r$ takes a
constant value.  Furthermore, $r$ carries $N_\F$ surjectively to
$W_\F$. \endproclaim

\demo{Proof} As in~(3.5), let $C_1$ be the canonical cross-section
$HK\cdot 1\cdot M_P U_P$.  Let $f_1$ be a marked lattice
corresponding to a point in $C_1$.  For this $f_1$, let $\frak A =
\frak A_{\F, f_1}$, and let $\bold t_1 \in
\frak A$ be the point constructed in Lemma~7.3.  We define
$$
N_\F = \{ \pmb\rho_1 \cdot f_1 \mid \pmb\rho_1 < \bold t_1 \},
\tag 7.4.1
$$
where $f_1$ ranges over all marked lattices representing points of
$C_1$, and $\bold t_1$ depends on $f_1$.  For any given $f_1$, the map $r$ is
constant on the set $\{ \pmb\rho_1 \cdot f_1 \} = N_\F \cap
\text{(fiber)}$, by Corollary~7.3. 

To prove $N_\F$ is an open set in $X$, we must show the $t_j$ vary
continuously with $f_1$.  The $\mu_i$ of~(4.1) vary continuously as a
function of $f$, by \cite{A3, \S3(i)}.  Hence the $\beta_j$ of~(7.3)
vary continuously with~$f$.  For any continuous family of lattices
(not necessarily normalized by $m(L) = 1$), the function $m(L)$ is
also continuous in $L$; applying this to $L^\dag$ (defined in~(7.3)),
we see the $\alpha_j$ of~(7.3) are continuous in~$f$.  Hence the $t_j$
are continuous in $f$.

By Corollary~7.2, $r$ carries $N_\F$ into $W_\F$.  In the rest of the
proof, we prove the map is surjective.  Let $f$ be a marked lattice
representing a point in $W_\F$.  Let $d_j = \dim_k V_j$.  Let $f_1$ be
the marked lattice corresponding to the point in $C_1$ in the
$A_P$-fiber through $f$; let $\pmb\rho \cdot f = f_1$.  Choose $\bold
t$ for $f_1$ as in Lemma~7.3.  For a sufficiently small $\tau \in
(0,1)$, $\pmb\tau = (\tau,\dots,\tau)$ satisfies $\pmb\tau < \bold t
\cdot \pmb\rho$, implying $\pmb\tau\cdot f \in N_\F$.  We will show
$r$ carries $\pmb\tau\cdot f$ back to $f$.

The lattice for $\pmb\tau\cdot f$ is the lattice for $f$ scaled by a
factor of $1/\tau^{j-1}$ on $\check V_j$.
The minimal vectors of $f$ span their intersection with $f(V_1)$, by
the definition of $W_\F$.  Since $1/\tau > 1$, the minimal vectors of
$\pmb\tau\cdot f$ are exactly the minimal vectors of $f$ that lie in
$f(V_1)$.  Hence the minimal vectors of $\pmb\tau\cdot f$ lie in
$f(V_1)$ and span $f(V_1)$.  This means the steps $\rs 1, \dots,
\rs{d_1 - 1}$ of the well-rounded retraction are trivial.

The retraction $\rs {d_1}$ fixes $f(V_1)$ and acts by a scalar
$\mu_{d_1} \le 1$ on all orthogonal directions.  Consider what happens
as a scalar $\mu$ descends from 1 to $\mu_{d_1}$.  The marked lattice
$\pmb\tau \cdot f$ becomes
$$
(\mu^{-1}\tau, \tau, \dots, \tau) \cdot f.
\tag 7.4.2
$$
As long as $1/(\mu^{-1} \tau) > 1$, there are no new minimal vectors.
But when $\mu^{-1} \tau = 1$, the minimal vectors of~(7.4.2) are
exactly the minimal vectors of $f$ that lie in $f(V_2)$.  These span
$f(V_2)$ (again by the definition of $W_\F$), which has rank bigger
than the rank of $f(V_1)$.  So by~(4.1.1) for $i=d_1$, we must have
$\mu_{d_1} = \tau$, and at the end of the $\rs {d_1}$ stage we have
produced $(1, \tau, \dots, \tau) \cdot f$.  The well-rounded
retraction goes on similarly to produce $(1, 1, \tau, \dots, \tau)
\cdot f$, \dots, and finally $(1,\dots,1) \cdot f = f$.  \qed\enddemo

\subheading{(7.5)} The neighborhood $N_\F$ defined in Proposition~7.4
lies in $X$, the interior of $\bar X$.  We now extend it to a
neighborhood $\bar N_\F$ of $e(P)$ in the corner $X(P)$.  As in~(3.7),
the identification $\frak A = \frak A_{\F,f_1} \cong A_P$ extends to
an identification $\bar \frak A = \bar \frak A_{\F, f_1} \cong \bar
A_P$, where $\bar \frak A$ is obtained from $\frak A$ by allowing one
or more of the $\rho_j$ in $(\rho_1, \dots, \rho_{j-1})\in \frak A$ to
become~0.  With the notation of~(7.4.1),

\definition{Definition} $\bar N_\F = \{ f = \pmb\rho_1 \cdot f_1 \mid
\pmb\rho_1 \in \bar\frak A_{\F, f_1}, \pmb\rho_1 < \bold t_1 \}$,
where $f_1$ ranges over all marked lattices representing points of
$C_1$, and $\bold t_1$ depends on $f_1$ as in Lemma~7.3.
\enddefinition

It is clear that $\bar N_\F$ is a tubular neighborhood of $e(P)$ in
$X(P)$, in the sense that $\bar N_\F$ is homeomorphic to $e(P) \times
[0,1)^{l-1}$ with $e(P)$ included as $e(P) \times (0,\dots,0)$.  Also,
the set $\bar N_\F$ is $(\Gamma_0\cap P)$-invariant.

\subheading{(7.6)} We now define a map $\bar r : \bar X \to W$ that
continuously extends the well-rounded retraction~$r$.  Choose
$x\in \bar X$.  If $x\in X$, let $\bar r(x) = r(x)$.  Otherwise, $x\in
e(P)$ for a unique $P$ with associated $\F$.
Hence $x\in \bar N_\F$.  Consider the fiber $\phi'$ of $\bar
N_\F$ over $x$.  By Proposition~7.4, $r$ takes a constant value on
$\phi' \cap X$.  Define $\bar r(x)$ to be this common value
$r(\phi' \cap X)$.

\proclaim{Proposition} The map $\bar r$ is a continuous
$\Gamma$-equivariant extension of $r$ to $\bar X$.  For any $\Q$-flag
$\F$, it is constant on the $\bar A_P$-fibers of $\bar N_\F$.  It is a
retraction $\bar X \to W$ (i.e.~it fixes $W$ pointwise). \endproclaim

\demo{Proof} It suffices to prove $\bar r$ is constant on $\phi'$.
For then, the continuity of $\bar r$ on $\bar N_\F$ will follow from
the tubular neighborhood structure and the continuity of $r$ on
$N_\F$.  The continuity of $\bar r$ on $\bar X$ will follow from its
continuity on $X$ and on each $\bar N_\F$, because these form an open
cover of $\bar X$.  The $\Gamma$-equivariance of $\bar r$ will follow
from that of~$r$.

Let $y$ lie on the fiber $\phi'$ of $\bar N_\F$ over $x$, with
$y\ne x$.  If $y \in X$, then $\bar r(x) = r(x) = \bar r(y)$ by
construction.  Next, suppose $y\in e(P_1)$ for some $P_1$
with associated $\F_1$.  Then $P_1 \supseteq P$.  Now
$A_{P_1}$ is naturally
a subgroup of $A_P$ by~(3.3), and $\bar A_{P_1} \subset \bar
A_P$.  So $\phi'$ fibers over $\bar e(P_1)$ with fibers
given by the action of $\bar A_{P_1}$.  Let $\phi_1'$ be the
intersection of $\phi'$ with the $\bar A_{P_1}$-fiber of
$\bar N_{\F_1}$ over $y$.  Then
$$
\align
\bar r(y) &= r(\phi_1' \cap X)
               \qquad\text{by def.~of $\bar r(y)$} \\
          &= r(\phi' \cap X)
               \qquad\foldedtext{since $\phi_1' \cap X
                                 \subseteq \phi' \cap
                                 X$, and $r$ is constant on the latter} \\
          &= \bar r(x) \qquad\text{by def.~of $\bar r(x)$}.
\endalign
$$
This proves $\bar r$ is constant on the $\bar A_P$-fibers of $\bar
N_\F$.

Finally, $\bar r$ fixes $W$ pointwise because $r$ does. \qed\enddemo


\heading Section 8---Extending the Deformation Retraction to the
Borel-Serre Boundary \endheading

\subheading{(8.1)} We recall from~\cite{S} some facts about Saper's
tiling of $\bar X$.  Saper constructs a $\Gamma$-equivariant subspace
$X_0$ of $X$ called the {\sl central tile.} This is a codimension-zero
submanifold with corners in $X$.  The closed boundary faces
$\partial^P X_0$ are in one-to-one correspondence with the proper
parabolic $\Q$-subgroups $P$; inclusions $P \subseteq
P'$ correspond to inclusions $\partial^P X_0 \subseteq
\partial^{P'} X_0$, and conversely.  Every $y\in \partial X_0$ lies on
$\partial^P X_0$ for a unique smallest $P$.  There is a
$\Gamma$-equivariant piecewise-analytic deformation retraction $r_S :
\bar X \to X_0$; more precisely, $r_S$ fits into a piecewise-analytic
family $r_{t, S}$ for $t\in [0,1]$ with $r_{0,S} = $(identity on $\bar
X$) and $r_{1,S} = r_S$.  The retraction $r_S$ is uniquely
characterized by
$$
r_S(y \geod a) = y \qquad\text{whenever }
    y\in \partial^P X_0 \text{ and } a \in \bar A_P, a\le 1.
\tag 8.1.1
$$
The converse of~(8.1.1)
also holds: if $y\in \partial^P X_0$, and $P$ is the smallest with
this property, then
$$
r_S(y') = y 
\quad\text{implies}\quad
\exists a \in \bar A_P, a \le 1, \text{ such that } y' = y \geod a.
\tag 8.1.2
$$

For any arithmetic subgroup $\Gamma$, the quotient $X_0/\Gamma$ is compact.
Also, for any $P$ and given any open neighborhood $U$ of $\bar
e(P)/(\Gamma\cap P)$ in $\bar X/\Gamma$, we may choose $X_0$ so that
$$
\partial^P X_0 / (\Gamma\cap P) \subseteq U. \tag 8.1.3
$$
In fact, given any
central tile, we may choose enlargements of it whose complements in
$\bar X$, modulo~$\Gamma$, form a fundamental system of neighborhoods
of $\partial \bar X/\Gamma$ in $\bar X/\Gamma$.

\subheading{(8.2)} We will now fix a certain $\Gamma$-invariant open
neighborhood $B$ of $\partial \bar X$ in $\bar X$.
There are only finitely many classes of $\Q$-flags mod $\Gamma$.  Let
$\F$ represent such a class.  Let $\bar N_\F$ be the tubular
neighborhood of $e(P)$ in $X(P)$ as in~(7.6), and let
$$
B = \bigcup_\F
       \bigcup_{\gamma\in \Gamma} \gamma\cdot \bar N_\F.
$$
Clearly $B/\Gamma$ is an open neighborhood of $\partial \bar X/\Gamma$ in
$\bar X/\Gamma$.

From now on, we fix a central tile $X_0$ with the following properties:
\roster
\item"(i)" $\partial^P X_0 \subseteq \bigcup_{\F'} \bigcup_{\gamma\in
\Gamma\cap P} (\gamma\cdot \bar N_{\F'})$ for any proper parabolic
$\Q$-subgroup $P$ with associated $\F$, where $\F'$ runs through a
set of $(\Gamma\cap P)$-representatives of $\Q$-flags containing~$\F$
(the number of such representatives is finite);
\item"(ii)" $\partial X_0 \subseteq B$, and $X_0 \cup B
= \bar X$;
\item"(iii)" the well-rounded retract $W$ lies in the interior of
$X_0$.
\endroster
By~(8.1.3), it is possible to arrange for~(i).  Part~(ii) follows
from~(i), since there are only finitely many $\F$ mod~$\Gamma$.  To
arrange for~(iii), recall that $W/\Gamma$ is compact, so that there is
a $\Gamma$-invariant
open neighborhood of $\partial \bar X$ in $\bar X$ whose
closure does not meet $W$; enlarge $X_0$ if necessary
so that its boundary lies in this neighborhood.

\subheading{(8.3)} We now construct the map $R_t$ described in~(0.3).

\definition{Definition} The map $R_t : \bar X \times [0,1] \to \bar X$
is given as follows.  If $t\in [0,\frac12]$, let $R_t(x) = r_{2t,
S}(x)$.  If $t\in [\frac12, 1]$, let $R_t(x) = r_{2t-1}(r_S(x))$.
\enddefinition

The map $R_t$ is well-defined after $t= \frac12$ because the image of
$r_S$ lies in $X$ (not $\partial \bar X$), and $r$ is defined on $X$.

\proclaim{Proposition} The map $R_t$ is a $\Gamma$-equivariant
deformation retraction of $\bar X$ onto $W$.  \endproclaim

\demo{Proof} $R_t$ is $\Gamma$-equivariant, since both $r_{t,S}$ and
$r_t$ have this property.  Since $W\subseteq X_0$ by~(iii) in~(8.2),
$W$ is fixed pointwise by $r_S$.  We know $W$ is fixed pointwise by
$r$.  Hence the image of $R_1$ is at least $W$.  On the other hand,
the image is contained in $W$ because the image of $r$ is $W$.
\qed\enddemo

\subheading{(8.4)} We now show that our global retraction $\bar r$
of~(7.6) is
really a deformation retraction.

\proclaim{Theorem} $\bar r = R_1$. \endproclaim

\demo{Proof} If $x\in X_0$, then $r_S$ does not move $x$, and $R_1(x)
= r(x) = \bar r(x)$.  Next, assume $x\in \bar X$ with $x \notin X_0$.
The retraction $r_S$ carries $x$ to a point $y$, and there is a unique
smallest $P$ such that $y$ lies on the face $\partial^P X_0$.
However, by~(8.1.2), $x$ is of the form $x = y \geod a$ for some $a\in
\bar A_P$, $a\le 1$. By~(8.2) part~(i), $y \in \bar N_{\F''}$ for some
$\F'' \supseteq \F$.  Since $\bar A_P \subseteq \bar A_{P''}$, we have
$x\in \bar N_{\F''}$, and $x$ lies on the same $\bar A_{P''}$-fiber of
$\bar N_{\F''}$ as $y$.  By~(7.6) applied to~$\cal F''$, we have $\bar
r(x) = \bar r(y)$.  But $\bar r(y) = r(y) = R_1(x)$ by the definition
of $R_t$.  \qed\enddemo

\subheading{(8.5)} Here is a simple consequence of Theorem~8.4.

\proclaim{Corollary} The map $\bar r : \bar X \to W$ is a $\Gamma$-equivariant
homotopy equivalence.  $\bar X/\Gamma$ and $W/\Gamma$ have the same
homotopy type via the map $\bar r$ mod~$\Gamma$.  \endproclaim

We now generalize this corollary to all the boundary components.

\proclaim{Theorem} Let $P$ be a proper parabolic $\Q$-subgroup, with
associated $\F$.  Then $\bar r$ induces a $(\Gamma\cap P)$-equivariant
homotopy equivalence of $\bar e(P)$ and $W_\F$. \endproclaim

\demo{Proof} We must define a $(\Gamma\cap P)$-equivariant homotopy
inverse $\bar s$ for $\bar r
\mid_{\bar e(P)}$.  Choose a point $x\in W_\F$, represented by a
marked 
lattice $f$.  Let $\pmb\tau = (\tau, \dots, \tau) \in \frak A = \frak
A_{\F,f}$.  We
define $\bar s(x)$ to be the point of $\bar X$ corresponding to
$$
\lim_{\tau\to 0+} \pmb\tau \cdot f.
\tag 8.5.1
$$
For all $\tau > 0$,  $\pmb\tau \cdot f$
lies in a single
$A_P$-fiber of $X$.  Let $q_P$ be as in~(3.5).  It
is clear that the limit in~(8.5.1) is
$q_P(x)$.  Thus an equivalent definition is
$$
\bar s \text{ is the restriction of } q_P \text{ to } W_\F.
\tag 8.5.2
$$
By~(8.5.2), $\bar s$ is continuous.  Since $q_P$ commutes with the action of
$P$, $\bar s$ is $(\Gamma\cap P)$-equivariant.

We now show that $\bar r \circ \bar s$ is the identity on $W_\F$.  For
a sufficiently small $\tau_0 > 0$, $(\tau_0, \dots, \tau_0) \cdot f$
will lie in $N_\F$.  By~(7.6), the value of $\bar r$ on the point
$(\tau, \dots, \tau) \cdot f$ is the same for all $\tau \in [0,
\tau_0]$.  Using~(8.5.1), $\bar r \circ \bar s(x) = \bar r((\tau_0,
\dots, \tau_0) \cdot f)$.  However, the proof of Proposition~7.4 shows
$\bar r ((\tau_0, \dots, \tau_0)\cdot f) = f$.

Next, we show $\bar s \circ \bar r \mid_{\bar e(P)}$ is homotopic to
the identity on $\bar e(P)$.  Consider the family of maps
$$
q_P \circ R_t \qquad\text{ restricted to $\bar e(P)$}.  \tag 8.5.3
$$
Here $R_t$ is from~(8.3), and $q_P$ is extended continuously to a
neighborhood of $\bar e(P)$ in $\bar X$ in a natural way \cite{B-S,
(5.5)}.  (A dangerous mistake would be to assume that the canonical
cross-sections in $X(P)$ extend over $\bar e(P)$ as cross-sections of
$q_P$.) The function in~(8.5.3) is continuous because the image of
$\bar e(P)$ under $R_t$ lies in $\bar e(P) \cup X$, where $q_P$ is
continuous.  At $t=0$, (8.5.3) is the identity $\bar e(P) \to \bar
e(P)$.  At $t=1$, for any $x\in \bar e(P)$,
$$
\align
q_P \circ R_1(x) &= q_P (\bar r(x)) \qquad\text{by (8.4)} \\
  &= \bar s(\bar r(x)) \qquad\text{by (8.5.2)}.
\endalign
$$
Thus $q_P \circ R_t$ affords a $(\Gamma\cap P)$-equivariant
homotopy between $\bar s \circ \bar r |_{\bar e(P)}$ and the
identity.  \qed\enddemo

\remark{Remark} Formula (8.5.3) defines a well-rounded retraction within
$\bar e(P)$ itself, carrying $\bar e(P)$ onto a subset of
$e(P)$ homeomorphic to $W_\F$.  We know that $\bar e(P) /
(\Gamma\cap P)$ is a compactification of a bundle whose base is a
locally symmetric space (perhaps with finite quotient singularities)
for a group of lower rank, and whose fiber is an arithmetic quotient
of $U_P$.  Accordingly, there is a cellular quotient map from $W_\F /
(\Gamma\cap P)$ to the well-rounded retract in the lower-rank locally
symmetric space.  This map may be constructed by dualizing the
techniques of~\cite{M2}; see~(10.7).
\endremark


\heading Section 9---The Main Results \endheading

In (9.5) we establish our main result, as we construct the spectral
sequence described in the introduction.  The bulk of the proofs
resides in (9.2)--(9.4).  We set up two spectral sequences,
one~(9.2.1) computing the cohomology of the Borel-Serre boundary
$\partial\bar X/\Gamma$, and the other~(9.3.1) based on the
subcomplexes $W_\F/(\Gamma\cap P)$ in the well-rounded retract
$W/\Gamma$.  Because the map $\bar r$ gives homotopy equivalences on
$\bar X/\Gamma$ and all its faces, it induces an isomorphism between
these spectral sequences, and between their abutments.  We set up
canonical maps $H^*(X/\Gamma) \to$~(9.2.1) and $H^*(W/\Gamma)
\to$~(9.3.1).  Finally, we show that the map $H^*(W/\Gamma)
\to$~(9.3.1) computes the canonical restriction map $H^*(\bar
X/\Gamma) \to H^*(\partial\bar X/\Gamma)$.

\subheading{(9.1)} Let $A^{p,q}$ be any first-quadrant double complex.
There are two standard filtrations of the double complex, $\{ A^{p,q}
\mid q \ge q_0\}$ and $\{A^{p,q} \mid p \ge p_0\}$.  The spectral
sequences these give are called respectively the {\sl type~I\/} and
{\sl type~II sequences for\/} $A^{p,q}$.  Our main results involve
type~II sequences.  The cohomology of the single complex associated to
$A^{p,q}$ is called the {\sl abutment\/} of (either) spectral
sequence.

Throughout the paper, (co)homology groups have coefficients in any
fixed abelian group.

\subheading{(9.2)} We set up a Mayer-Vietoris spectral sequence for
the cohomology of $\partial\bar X/\Gamma$.
Let $\Phi_l$ be a set of representatives of the $\Gamma$-equivalence
classes of the $\Q$-flags $\cal F = \{0 \subsetneqq V_1 \subsetneqq
\cdots \subsetneqq V_l = k^n \}$ (see~(3.2)).  There are only finitely
many such equivalence classes.  Let $\cal E = \{ \bar e(P)/(\Gamma\cap P) \mid \F \in \Phi_2 \}$ (corresponding to maximal
proper parabolic $\Q$-subgroups).  This $\cal E$ is a cover of
$\partial\bar X/\Gamma$.  By (3.5), the non-empty $(p+1)$-fold
intersections of distinct members of $\cal E$ are exactly the sets in
$\{ \bar e(P)/(\Gamma\cap P) \mid \F \in \Phi_{p+2} \}$.  This
means we are in the correct setting to use the \v Cech cohomology
techniques of \cite{B-T, \S15}.

For $0\le p \le n-2$, define
$$
\cal X^{p,q} =
  \bigoplus_{\F \in \Phi_{p+2}} C^q(\bar e(P) / (\Gamma \cap P)),
$$
where $C^q$ denotes the singular $q$-cochains.  We make $\cal X^{p,q}$
into a double complex, where the vertical maps are $(-1)^p$ times the
coboundary maps, and the horizontal maps are induced as in \v Cech
cohomology from the inclusions $\bar e(P) \hookleftarrow \bar
e(P')$ 
for $\F' \supseteq \F$, with the alternating sign convention of
\cite{B-T, (15.7.1)}.  Exactly as in \cite{B-T, (15.7)}, we see that
the $q$-th row of $\cal X^{p,q}$ is exact, except at the $p=0$
position where the kernel equals $C^q(\partial\bar X/\Gamma)$.  Thus the
type~I spectral sequence for $\cal X^{p,q}$ collapses at $E_2$ to
$H^*(\partial\bar X/\Gamma)$.  Therefore, the total complex of $\cal
X^{p,q}$ computes $H^*(\partial\bar X/\Gamma)$.

Let $E^{p,q}_{r,\cal X}$ be the type~II spectral sequence for
$\cal X^{p,q}$.  We have
$$
E^{p,q}_{1,\cal X}
  = \bigoplus_{\F \in \Phi_{p+2}} H^q(\bar e(P) / (\Gamma \cap P))
  \Rightarrow  H^{p+q}(\partial\bar X/\Gamma).
\tag 9.2.1
$$

\subheading{(9.3)} For $0\le p \le n-2$, define
$$
\cal W^{p,q} =
  \bigoplus_{\F \in \Phi_{p+2}} C^q(W_\F / (\Gamma \cap P)),
$$
made into a double complex in the same way as $\cal X^{p,q}$.
(Recall that $W_{\F} \hookleftarrow W_{\F'}$ whenever $\F' \supseteq
\F$.)    Let $E^{p,q}_{r,\cal W}$ be the type~II spectral
sequence for $\cal W^{p,q}$.  As above,
$$
E^{p,q}_{1,\cal W}
  = \bigoplus_{\F \in \Phi_{p+2}} H^q(W_\F / (\Gamma \cap P)).
\tag 9.3.1
$$

Let $\cal W^*$ denote the single complex associated to $\cal W^{p,q}$.
Thus the abutment of~(9.3.1) is $H^{p+q}(\cal W^*)$.

\subheading{(9.4)} We now define a map of double complexes.  Consider
the single-column double complex
$$
\hat \cal X^{p,q} =
\left\{
\matrix
\format\c&\qquad\l \\
C^q(\bar X/\Gamma) &\text{if }p = 0 \\
0 &\text{if }p > 0.
\endmatrix
\right.
$$
There is an obvious map $\hat \cal X^{p,q} \to \cal
X^{p,q}$, where for $p = 0$ we map $\omega \in C^q(\bar X/\Gamma)$ to the
direct sum over $\F\in \Phi_2$ of the restriction of $\omega$ to
$C^q(\bar e(P) / (\Gamma\cap P))$.

\proclaim{Lemma} The map $\hat \cal X^{p,q} \to \cal X^{p,q}$ is a
chain map of double complexes. \endproclaim

\demo{Proof} Let $\omega \in C^q(\bar X/\Gamma)$, and let $\tau$ be
its image in $\cal X^{0,q}$.  The only thing to check is that the
horizontal arrow $\delta : \cal X^{0,q} \to \cal X^{1,q}$ carries
$\tau$ to~0.  But as we have said, the kernel of $\delta$ is precisely
the set of elements that come from a global cochain on $\partial \bar
X/\Gamma$.  And $\tau$ does come from a global cochain, namely
the restriction of $\omega$ to $\partial \bar X/\Gamma$.  \qed\enddemo

Considering the type~I spectral sequences, we see that $\hat \cal
X^{p,q} \to \cal X^{p,q}$ induces on the abutments the canonical
restriction map $H^q(\bar X/\Gamma) \to H^q(\partial \bar X/\Gamma)$.

Next, consider
$$
\hat \cal W^{p,q} =
\left\{
\matrix
\format \c & \qquad\l \\
C^q(W/\Gamma) &\text{if }p = 0 \\
0 &\text{if }p > 0,
\endmatrix
\right.
$$
made into a double complex the same way $\hat \cal X$ was.
There is an obvious map $\psi : \hat \cal W^{p,q} \to \cal W^{p,q}$ carrying
$\omega \in C^q(W/\Gamma)$ to the direct sum over $\F\in \Phi_2$ of the
restrictions of $\omega$ to $C^q(W_\F / (\Gamma\cap P))$.  Again,
$\psi: \hat
\cal W^{p,q} \to \cal W^{p,q}$ is a chain map of double complexes.

\subheading{(9.5)} We now give the main result of the paper.

\proclaim{Theorem} (i) We have a commutative diagram
$$
\CD
H^*(\bar X/\Gamma) @>\text{restriction}>> H^*(\partial\bar X/\Gamma) \\
@V{\cong}VV                               @VV{\cong}V \\
H^*(W/\Gamma) @>{\psi^*}>>                    H^*(\cal W^*)
\endCD
$$
where the vertical maps are natural isomorphisms and the top map is
the canonical restriction map.

(ii) There is a spectral sequence
$$
E_{1,\cal W}^{p,q} = \bigoplus_{\F \in \Phi_{p+2}} H^q(W_\F / (\Gamma\cap P))
     \Rightarrow H^{p+q}(\cal W^*),
\tag 9.5.1
$$
and $\psi$ induces a map of spectral sequences given on the $E_1$ page
by the natural restriction map $H^q(W/\Gamma) \to E_{1,\cal W}^{0,q}$.
\endproclaim

\noindent Here $X$ is as in~(1.5), $\Gamma$ is any arithmetic subgroup
of $\Gamma_0$ as in~(2.2), $\Phi_l$ is as in~(9.2), and $P$ is the
$\R$-points of the parabolic $\Q$-subgroup associated to~$\F$ as
in~(3.2).  The space $W_\F$, a closed subcomplex of the well-rounded
retract $W$, is defined in~(6.1).  We use any fixed abelian group of
coefficients.

\demo{Proof} First we observe that the global deformation retraction
$\bar r: \bar X \to W$ induces a map of double
complexes $\cal X^{p,q} \to \cal W^{p,q}$.  The map $\bar r$ commutes
with the vertical maps of these double complexes, because any
continuous map commutes with the coboundary maps in singular
cohomology.  More importantly, $\bar r$ commutes with the horizontal
maps because it is defined {\it globally\/} on $\bar X$, and because the
horizontal maps are induced from inclusion maps on various subspaces
of~$\bar X$. 

It follows that $\bar r$ induces a homomorphism of spectral sequences
$E^{p,q}_{r, \cal X} \to E^{p,q}_{r, \cal W}$.  The map on $E_1$ terms
comes from the maps $H^q(\bar e(P) / (\Gamma\cap P)) \to H^q(W_\F
/(\Gamma\cap P))$ induced for all $\F$ by $\bar r$.  By
Theorem~8.5, the latter maps are canonical isomorphisms.  By
\cite{Br, VII.2.6}, the right-hand arrow in~(i) is an isomorphism.

Similarly, $\bar r$ induces a map $\hat \cal X^{p,q} \to \hat \cal
W^{p,q}$.  By Corollary~8.4, it gives a canonical isomorphism in
cohomology.  This is the left vertical arrow in~(i).

Part~(9.4) defined the maps of double complexes that induce the
horizontal arrows in~(i), and it explained why the top arrow is the
canonical restriction.  This proves~(i).

The spectral sequence in~(ii) is~(9.3.1).  The assertion about~$\psi$
in~(ii) is immediate from the definition of~$\psi$.  \qed\enddemo


\heading Section 10---Practical Comments and Related Results
\endheading

\subheading{(10.1)} We make some remarks about how to use the
spectral sequence~(9.5.1) in practice.  When we work with the
finite cell complex $W/\Gamma$, we represent each cell by some piece
of data.  (See~(10.6), \cite{M-M0}, or \cite{A-G-G}.) Since
$W_{\cal F}$ is a subcomplex of $W$, it is tempting to imagine
representing each cell of $W_{\cal F}/ (\Gamma \cap P)$ by some piece
of the $W/\Gamma$ data.  However, two cells of $W_{\cal F}$ may be
equivalent mod~$\Gamma$ and yet inequivalent mod $(\Gamma \cap P)$.
An example when $\Gamma = \GL_3(\Z)$, $\F = \{ 0 \subseteq (*,*,0)
\subseteq \Q^3 \}$, is provided by the two 0-cells of $W$ whose minimal
vectors are given by the columns of the following matrices (together
with their negatives):
$$
\pmatrix
1 & 0 & 1 & 0 & 0 & -1 \\
0 & 1 & 1 & 0 & 1 & 0 \\
0 & 0 & 0 & 1 & 1 & 1
\endpmatrix
\quad\text{ and }\quad
\pmatrix
1 & 0 & 0 & 1 & 1 & 0 \\
0 & 1 & 0 & 0 & 1 & 1 \\
0 & 0 & 1 & 1 & 1 & 1
\endpmatrix.
$$
The cells are $\Gamma$-equivalent.  (In the language of \cite{M-M0},
their $\cal C$-configurations are both complete quadrilaterals.)  But
any element of $\Gamma$ preserving
$\F$ must carry the left-hand matrix to a matrix with three 0's in its
bottom row.

In general, each $\Gamma$-class of cells of $W_{\cal F}$ will break up
into finitely many $(\Gamma \cap P)$-classes.  So a cell in $W_{\cal
F}/ (\Gamma \cap P)$ is represented by a piece of data for $W/\Gamma$,
plus a finite amount of extra structure depending on $\cal F$.  Say
one has worked out all the classes in~(9.5.1) as explicit cocycles on
the cell complexes $W_{\cal F}/ (\Gamma \cap P)$.  To compute the
$d_1$ and higher differentials, one has to work out how the
data-with-extra-structure in $W_{\cal F}/ (\Gamma \cap P)$ is
connected with the data-with-extra-structure in $W_{\cal F'}/ (\Gamma
\cap P)$ for $\cal F' \supseteq \F$.  {\smc Sheafhom}, a suite of
object-oriented programs currently under development by the second
author, may help to automate such computations.

Say that $\Gamma$ is {\sl small enough\/} if, whenever $\cal F$, $\cal
F'$ are distinct $\Q$-flags and $W_{\cal F} \cap W_{\cal
F'} \ne \varnothing$, there is no element of $\Gamma$ carrying $\cal
F$ to $\cal F'$.  When $\Gamma$ is small enough in this sense, we need
never worry about the extra structure just described.  To see why,
assume there were two cells $C_a, C_b \in W_{\cal F}$ which were not
$(\Gamma \cap P)$-equivalent but such that $C_b = \gamma \cdot C_a$
for some $\gamma \in \Gamma$.  Set $\cal F' = \gamma^{-1} \cal F$; we
have $\F' \ne \F$ since $\gamma \notin P$.  Because $C_b$ respects
$\cal F$, $C_a$ must respect $\cal F'$.  Hence $C_a \in W_{\cal F}
\cap W_{\cal F'}$.  Since $\Gamma$ was assumed small enough, the
existence of such a $\gamma$ is a contradiction.

When $k=\Q$ and $n\le 4$, one can use the methods of \cite{M1} to show
that the principal congruence subgroups $\Gamma(N) \subset \GL_n(\Z)$
of prime level $N\ge3$ are small enough in the above sense.  We do not
know whether all neat $\Gamma$ are small enough.

\subheading{(10.2)} For convenience, we record the homology version
of our main theorem.  We use the obvious analogues of the notation
of~(9.5).

\proclaim{Theorem} (i) We have a commutative diagram
$$
\CD
H_*(\partial\bar X/\Gamma) @>\text{inclusion}>> H_*(\bar X/\Gamma) \\
@V{\cong}VV                               @VV{\cong}V \\
H_*(\cal W_*) @>{\psi_*}>>                    H_*(W/\Gamma) 
\endCD
$$
where the vertical maps are natural isomorphisms and the top map is
the canonical inclusion map.

(ii) There is a spectral sequence
$$
E^1_{p,q} = \bigoplus_{\F \in \Phi_{p+2}} H_q(W_\F / (\Gamma\cap P))
     \Rightarrow H_{p+q}(\cal W_*),
$$
and~$\psi$ induces a map of spectral sequences given
on the $E_1$ page by the natural inclusion map $E^1_{0,q} \to
H_q(W/\Gamma)$.
\endproclaim

This is proved by dualizing the constructions of Section~9.

\subheading{(10.3)} Consider one Borel-Serre boundary face $\bar
e(P)/(\Gamma\cap P)$.  As we mentioned in the introduction, our methods allow
us to compute the canonical maps $H^*(\bar X/\Gamma) \to H^*(\bar
e(P)/(\Gamma\cap P))$ and $H_*(\bar e(P)/(\Gamma\cap P)) \to H_*(\bar
X/\Gamma)$ without any spectral sequences.  We simply use the cellular
map $W_{\cal F}/(\Gamma \cap P) \to W/\Gamma$.

In the application to $\GL_4(\Z)$ described near the end of~(0.2), we
are primarily interested in the degree-five homology.
We believe that in degree five, all the homology of the boundary comes
from the maximal boundary faces.  We could therefore find
the image of
$$
\left( \bigoplus
         \Sb \cal F \text{ for maximal} \\ \text{faces mod }\Gamma \endSb
         H_5( W_{\cal F} / (\Gamma \cap P) ) \right)
  \to H_5(W/\Gamma)
$$
without computing any spectral sequences.

It would very interesting to find the image of the map $H_*(\partial
\bar X/\Gamma) \to H_*(\bar X/\Gamma)$ without having to find
$H_*(\partial \bar X/\Gamma)$ directly, and without having to work out
any spectral sequences.

\subheading{(10.4)} All our results hold for $\Gamma$-equivariant
cohomology.  That is, one may replace $H^*((\dots)/\Gamma)$ with
$H^*_\Gamma((\dots))$ in the theorems of Section~9 and throughout the
paper.  Of course, if we use coefficients in a field of
characteristic~0 or of characteristic prime to the order of any
torsion element of $\Gamma$, then the two kinds of cohomology are
canonically isomorphic.  All these results hold for homology.

In the $\Gamma$-equivariant setting, the summands in each column of
our spectral sequences would be replaced by $H^*_{(\Gamma\cap
P)}(W_\F)$, and each of these would have to be computed by a spectral
sequence in its own right.  We would get a spectral sequence of
spectral sequences.

Both $H^*_\Gamma(X)$ and $H^*(X/\Gamma)$ are interesting objects from
the point of view of number theory and automorphic forms.  Both have
Hecke operators acting on them.  In~\cite{A4}, the first author has
conjectured that any Hecke eigenclass in $H^*_\Gamma(X; \Bbb F_p)$,
for any~$p$, should have an attached mod~$p$ representation of the
absolute Galois group of~$\Q$.  In~\cite{A-M2}, we have conjectured
the same for the Hecke module $H^*(X/\Gamma; \Bbb F_p)$.  Both papers
provide examples.

\subheading{(10.5)} The method of this paper can easily be extended to
$\SL_n$.  If $G$ were $\SL_n(S)$, then $K\bs G$ would be the product
of the irreducible symmetric spaces for the $\SL_n(k_v)$.  In the case
$G = \GL_n(S)$ of this paper, $HK\bs G$ is the product of the
$\SL_n(k_v)$ symmetric spaces and a Euclidean factor $H \bs \left(
\prod_v \R_+ \right) \cong \R_+^{r_1+r_2-1}$, where $r_1, r_2$ are the
number of real and complex places of~$k$, respectively.

\subheading{(10.6)} Let $k=\Q$, so that $\Gamma \subseteq \GL_n(\Z)$
or $\SL_n(\Z)$, and let $n\le 4$.  In \cite{M1} and the survey article
\cite{M-M0}, an interpretation of the cell structure on the
well-rounded retract $W$ is given.  There is a one-to-one
correspondence in which each cell corresponds to a configuration of
points and lines (called a {\sl $\cal C$-configuration}) in the
projective space $\P^{n-1}(\Q)$.  It is easy to interpret our $W_{\cal
F}$ in this language.

\subheading{(10.7)} Again let $k=\Q$ and $n\le 4$.  Let $\breve
X/\Gamma$ be any Satake compactification of $X/\Gamma$.  In \cite{M2},
one constructs a cell structure on $\breve X/\Gamma$.  Intersecting
with $X/\Gamma$ gives a ``locally-compact cell structure'' on
$X/\Gamma$ that is dual to the well-rounded retract $W/\Gamma$ inside
$X/\Gamma$ (see \cite{M1, \S2} for definitions and precise
statements).  There is a cellular map from the barycentric subdivision
of the boundary of $\breve X/\Gamma$ to the barycentric subdivision of
$W/\Gamma$.  This offers an interesting contrast to the present paper.


\enddocument